\documentclass[ejs]{imsart}
\setattribute{journal}{name}{}
\RequirePackage[OT1]{fontenc}
\RequirePackage{amsthm,amsmath}
\RequirePackage{natbib}
\RequirePackage[colorlinks,citecolor=blue,urlcolor=blue]{hyperref}

\usepackage{amsfonts, amsopn, amssymb}
  \numberwithin{equation}{section}
\usepackage{color}
\usepackage{graphicx, subcaption, epstopdf}

\startlocaldefs
\numberwithin{equation}{section}
\theoremstyle{plain}
\newtheorem{theorem}{Theorem}[section]
\endlocaldefs

\newtheorem{assumption}[theorem]{Assumption}
\newtheorem{corollary}[theorem]{Corollary}
\newtheorem{definition}[theorem]{Definition}
\newtheorem{lemma}[theorem]{Lemma}

\theoremstyle{remark}

\newcommand{\BALD}{\begin{aligned}}
\newcommand{\EALD}{\end{aligned}}
\newcommand{\BALDS}{\begin{aligned*}}
\newcommand{\EALDS}{\end{aligned*}}
\newcommand{\BCAS}{\begin{cases}}
\newcommand{\ECAS}{\end{cases}}
\newcommand{\BEAS}{\begin{eqnarray*}}
\newcommand{\EEAS}{\end{eqnarray*}}
\newcommand{\BEQ}{\begin{equation}}
\newcommand{\EEQ}{\end{equation}}
\newcommand{\BIT}{\begin{itemize}}
\newcommand{\EIT}{\end{itemize}}
\newcommand{\BMAT}{\begin{bmatrix}}
\newcommand{\EMAT}{\end{bmatrix}}
\newcommand{\BNUM}{\begin{enumerate}}
\newcommand{\ENUM}{\end{enumerate}}
\newcommand{\eg}{e.g.}
\newcommand{\Eg}{E.g.}
\newcommand{\ie}{i.e.}

\newcommand{\BA}{\begin{array}}
\newcommand{\EA}{\end{array}}

\newcommand{\ones}{\mathbf 1}

\newcommand{\reals}{\mathbf{R}}

\DeclareMathOperator*{\argmin}{\arg\min}

\DeclareMathOperator{\Expect}{\mathbf{E}}

\DeclareMathOperator*{\minimize}{minimize}

\renewcommand{\Pr}{\mathbf{Pr}}
\DeclareMathOperator{\rank}{rank}

\DeclareMathOperator{\sign}{sign}

\DeclareMathOperator{\tr}{tr}

\newcommand{\pc}{\hspace{1pc}}

\newcommand{\abs}[1]{\left| #1 \right|}

\DeclareMathOperator{\linspan}{span}
\newcommand{\norm}[1]{\left\| #1 \right\|}

\newcommand{\iid}{\emph{i.i.d.}}

\DeclareMathOperator{\conv}{conv}

\DeclareMathOperator{\ext}{ext}

\DeclareMathOperator{\relint}{relint}

\newcommand{\bDelta}{\bar{\Delta}}

\newcommand{\brho}{\bar{\rho}}
\newcommand{\btheta}{\bar{\theta}}
\newcommand{\bvarrho}{\bar{\varrho}}

\DeclareMathOperator{\Fro}{F}
\newcommand{\hDelta}{\hat{\Delta}}
\newcommand{\hSigma}{\hat{\Sigma}}
\newcommand{\htheta}{\hat{\theta}}
\newcommand{\hTheta}{\hat{\Theta}}
\DeclareMathOperator{\IC}{IC}

\newcommand{\Thetas}{\Theta^\star}
\newcommand{\thetas}{\theta^\star}
\DeclareMathOperator{\op}{sp}

\newcommand{\bP}{\mathbf{P}}

\newcommand{\bU}{\bar{U}}
\newcommand{\bV}{\bar{V}}
\newcommand{\bz}{\bar{z}}
\newcommand{\cA}{\mathcal{A}}

\newcommand{\cG}{\mathcal{G}}

\newcommand{\cS}{\mathcal{S}}
\newcommand{\cX}{\mathcal{X}}
\newcommand{\hU}{\hat{U}}
\newcommand{\hV}{\hat{V}}
\newcommand{\hz}{\hat{z}}
\newcommand{\Mp}{M^\perp}
\newcommand{\tU}{\tilde{U}}
\newcommand{\tV}{\tilde{V}}

\begin{document}

\begin{frontmatter}

\title{On model selection consistency of regularized M-estimators}
\runtitle{Model selection consistency of regularized M-estimators}

\begin{aug}
\author{\fnms{Jason D.} \snm{Lee}\ead[label=e1]{jdl17@stanford.edu}},
\author{\fnms{Yuekai} \snm{Sun}\ead[label=e2]{yuekai@stanford.edu}},
\and
\author{\fnms{Jonathan E.} \snm{Taylor}\ead[label=e3]{jonathan.taylor@stanford.edu}}

\runauthor{Lee, Sun, and Taylor}

\affiliation{Stanford University}

\address{Institute for Computational and Mathematical Engineering\\
Stanford University \\
\printead{e1},
\printead*{e2}
}

\address{Department of Statistics\\
Stanford University \\
\printead{e3}
}
\end{aug}

\begin{abstract}
Regularized M-estimators are used in diverse areas of science and engineering to fit high-dimensional models with some low-dimensional structure. Usually the low-dimensional structure is encoded by the presence of the (unknown) parameters in some low-dimensional model subspace. In such settings, it is desirable for estimates of the model parameters to be \emph{model selection consistent}: the estimates also fall in the model subspace. We develop a general framework for establishing consistency and model selection consistency of regularized M-estimators and show how it applies to some special cases of interest in statistical learning. Our analysis identifies two key properties of regularized M-estimators, referred to as geometric decomposability and irrepresentability, that ensure the estimators are consistent and model selection consistent.
\end{abstract}



\end{frontmatter}

\section{Introduction}
\label{sec:intro}

The principle of parsimony is used in many areas of science and engineering to promote ``simple'' models over more complex ones. In machine learning, signal processing, and high-dimensional statistics, this principle motivates the use of sparsity inducing penalties for model selection and signal recovery from incomplete/noisy measurements. Usually the ``simplicity'' of the model is encoded by the presence of the (unknown) parameters in some low-dimensional model subspace, and it si desirable for estimates of the parameters to fall in the model subspace. This notion of correctness is termed \emph{model selection consistency.} In this work, we consider regularized M-estimators of the form
\BEQ
\minimize_{\theta\in E}\,\ell(\theta) + \lambda\rho(\theta),
\label{eq:regularized-m-estimator}
\EEQ
where $\ell$ is a convex, twice continuously differentiable loss, $\rho$ is a penalty function, and $E\subseteq\reals^p$ is a subspace. We identify two key properties of regularized M-estimators, referred to as geometric decomposability and irrepresentability, that ensure the estimators are consistent and model selection consistent. We also develop a general framework for analyzing the consistency and model selection consistency of M-estimators with geometrically decomposable penalties. When specialized to various statistical models, our framework yields some known and some new model selection consistency results.

The article is organized as follows: First, we review existing work on consistency and model selection consistency of regularized M-estimators. Then, in Section \ref{sec:geometrically-decomposable}, we describe geometrically decomposable penalties. Section \ref{sec:main-result} is devoted to our main result and some discussion of its consequences.  converse results on the necessity of the irrepresentable condition in Section \ref{sec:weakly-decomposable}. The final section, Section \ref{sec:examples}, is devoted to applications of our main result to various statistical models, including sparse regression and low-rank multivariate regression.

\subsection{Notation} Given a set $S\subset\reals^p$ and a point $x\in\reals^p,$ we use $P_S(x)$ to denote the \emph{projector} of $x$ on $\linspan(S):$
\[\textstyle
P_S(x) = \argmin_{y\in \linspan(S)}\frac12\norm{x - y}_2^2.
\]
Since $P_S(x)$ is a linear mapping, we write $P_Sx = P_S(x),$ where $P_S\in\reals^{p\times p}$ We use $B_q$ to denote the $q$ norm ball $\{x\in\reals^p\mid\norm{x}_q \le 1\}.$ For a semi-norm $\rho,$ we use $\rho^*$ to denote its \emph{dual semi-norm:}
\[\textstyle
\rho(x)^* = \sup_{\rho(x)\,\le\, 1} y^Tx.
\]
Finally, given a matrix $X\in\reals^{p_1\times p_1},$ we use $X^\dagger$ to denote its \emph{Moore-Penrose pseudoinverse.}

\subsection{Consistency of regularized M-estimators}

A large body of work in high-dimensional statistics focuses on obtaining sufficient conditions for consistency of regularized M-estimation. A recurring theme in this avenue of research is the notion of \emph{restricted strong convexity}. We refer to the past work section in \cite{negahban2012unified} and \cite{buhlmann2011statistics} for a comprehensive treatment of recent work on this topic.

\cite{negahban2012unified} proposes a unified framework for establishing consistency and convergence rates for M-estimators with penalties $\rho$ that are \emph{decomposable} with respect to a pair of subspaces $M,\bar{M}$:
$$
\rho(x + y) = \rho(x) + \rho(y),\textnormal{for all }x\in M,y\in\bar{M}^\perp.
$$
Many common penalties such as the lasso, group lasso, and nuclear norm are decomposable in this sense. \cite{negahban2012unified} also develop a general notion of restricted strong convexity and prove a general result that establishes the consistency of M-estimators with decomposable penalties. Using their framework, they establish consistency results for different statistical models including sparse and group sparse linear regression.

More recently, \cite{vandegeer2012weakly} proposes the notion of \emph{weakly decomposability}. A penalty $\rho$ is weakly decomposable if there is some norm $\rho_{\cS^c}$ on $\reals^{p - \abs{\cS}}$ such that $\rho$ is superior to the sum of $\rho$ and $\rho_{\cS^c}$; \ie
$$
\rho(x) \ge \rho(x_\cS) + \rho_{\cS^c}(x_{\cS^c}),\textnormal{for all }x\in\reals^p,
$$
where $\cS\subset[\,p\,]$ and $x_\cS\in\reals^{\abs{\cS}},x_{\cS^c}\in\reals^{p-\abs{\cS}}.$ Many common sparsity inducing penalties, including the $\ell_2/\ell_1$-norm (with possibly overlapping groups), are weakly decomposable. \cite{vandegeer2012weakly} shows oracle inequalities for the $\ell_1$ penalty generalizes to weakly decomposable penalties.

In the parallel world of signal processing, there is a rich literature on constrained M-estimators of the form
\BEQ
\minimize_{\theta\in\reals^p}\,\norm{\theta}_\cA\text{ subject to }\theta \in C,
\label{eq:constrained-m-estimator}
\EEQ
where $C\subset\reals^p$ is a convex set. \cite{candes2012simple} proposed a unified analysis of \eqref{eq:constrained-m-estimator} when $\norm{\cdot}_{\cA}$ is decomposable. By \cite{candes2012simple}, Definition 1, $\norm{\cdot}_{\cA}$ is decomposable at $\theta^\star\in\reals^p$ if $\partial\norm{\thetas}_\cA$ has the form
\[
\partial\norm{\thetas}_\cA = \{z\in\reals^p\mid P_T(z) = e,\rho(P_{T^\perp}(z))^* \le 1\}
\]
for some subspace $T\subset\reals^p$ and a point $e\in T.$ Above, $P_T$ (resp. $P_{T^\perp}$) is the orthogonal projection onto $T$ (resp. $T^\perp$). We emphasize the notion of decomposability by \cite{candes2012simple} is different from the notion by \cite{negahban2012unified}. \cite{candes2012simple} show exact recovery results for sparse/block-sparse vectors and low-rank matrices from random linear measurements depend upon the decomposability of the $\ell_1,\ell_2/\ell_1,$ and nuclear norms.

Recently, \cite{chandrasekaran2012convex} proposed the notion of an \emph{atomic norm}:
\[
\norm{x}_A = \inf\left\{t > 0\mid x\in t\conv(A)\right\}\text{ for a set of atoms }A.
\]
They develop a general framework for deriving both exact (in the noise-free case) and robust (in the noisy case) recovery results from random Gaussian measurements by solving convex optimization problems of the form \eqref{eq:constrained-m-estimator}.

The model selection consistency of regularized M-estimators has also been extensively studied. The most commonly studied problems are
\BNUM
\item sparse regression (including generalized linear models): \cite{zhao2006model,bunea2008honest,wainwright2009sharp,obozinski2011support,vaiter2011robust}
\item sparse covariance estimation and (more generally) structure learning: \cite{meinshausen2006high, kolar2010estimating,ravikumar2010high, jalali2011learning,loh2012structure}.
\ENUM
In addition to restricted strong convexity, these results also rely upon the notion of \emph{irrepresentability} originally proposed by \cite{zhao2006model}\footnote{An equivalent notion, called \emph{neighborhood stability}, was proposed by \cite{meinshausen2006high}.}. Despite extensive work on this area, there is no general framework for establishing model selection consistency of commonly used M-estimators.

\section{Geometrically decomposable penalties}
\label{sec:geometrically-decomposable}

Let $C\subset\reals^p$ be a closed convex set. Then the \emph{gauge function} on $C$ is
$$\textstyle
\gamma_C(x) = \inf_x\,\{\lambda\in\reals_+\mid x\in\lambda C\},
$$
and the \emph{support function} on $C$ is
\BEQ\textstyle
h_C(x) = \sup_y\,\{y^Tx\mid y\in C\}.
\label{eq:support-function}
\EEQ
Both gauge support functions are sublinear and should be thought of as semi-norms. If $C$ is a norm ball, \ie{} $C = \{x\mid\norm{x} \le 1\}$, then $\gamma_C$ is the norm and $h_C$ is the dual norm given by
\[\textstyle
\norm{y}_* = \sup_x\,\{x^Ty\mid\norm{x}\le 1\}.
\]

The support function is a supremum of linear functions, hence the subdifferential consists of the linear functions that attain the supremum:
\BEQ
\partial h_C(x) = \{y\in C\mid y^Tx = h_C(x)\}.
\label{eq:support-function-subdifferential}
\EEQ
The support function (as a function of the convex set $C$) is also additive over Minkowski sums, \ie{} if $C$ and $D$ are convex sets, then
$$
h_{C+D}(x) = h_C(x) + h_D(x).
$$
We use this property to express penalty functions as sums of support functions. \Eg\ if $\rho$ is a norm and the dual norm ball can be expressed as a (Minkowski) sum of convex sets $C_1,\dots,C_k$, then $\rho$ can be expressed as a sum of support functions:
$$
\rho(x) = h_{C_1}(x) + \dots + h_{C_k}(x).
$$

If a penalty $\rho$ can be expressed as
\BEQ
\rho(\theta) = h_A(\theta) + h_I(\theta) + h_{E^\perp}(\theta),
\label{eq:geometrically-decomposable-penalty}
\EEQ
where $A,I\subset\reals^p$ are closed convex sets and $E\subset\reals^p$ is a subspace, then we say $\rho$ is a \emph{geometrically decomposable} penalty. This form is general; if $\rho$ can be expressed as a sum of support functions, \ie{}
$$
\rho(\theta) = h_{C_1}(\theta) + \dots + h_{C_k}(\theta),
$$
then we can set $A$, $I$, and $E^\perp$ to be sums of the sets $C_1,\dots,C_k$ to express $\rho$ in geometrically decomposable form \eqref{eq:geometrically-decomposable-penalty}.
In many cases of interest, $A+I$ is a norm ball and $h_{A+I} = h_A + h_I$ is the dual norm. In our analysis, we further assume
\BNUM
\item $A,I$ are bounded.
\item $I$ contains a relative neighborhood of the origin, \ie{} $0\in\relint(I).$
\ENUM
To allow for unregularized parameters, we do not assume $A + I$ contains a neighborhood of the origin. Thus $\rho$ is not necessarily a norm. We summarize the form of geometrically decomposable penalties in a definition.

\begin{definition}
\label{def:geometrically-decomposable}
A regularizer is \emph{geometrically decomposable} in terms of convex sets $A,I\subset\reals^p$ and a subspace $E\subset\reals^p$ if
\[
\rho(\theta) = h_A(\theta) + h_I(\theta) + h_{E^\perp}(\theta).
\]
We assume $A,I$ are bounded and $0\in\relint(I).$
\end{definition}

The notation $A,I$ should be as read as ``active'' and ``inactive'': $\linspan(A)$ should contain the (unknown) parameter vector and $\linspan(I)$ should contain deviations that we want to penalize. For example, if we know the sparsity pattern of the unknown parameter vector, then $A$ should span the subspace of all vectors with the correct sparsity pattern. 

The third term enforces a subspace constraint $\theta\in E$ because the support function of a subspace is the (convex) indicator function of the orthogonal complement:
$$
h_{E^\perp}(\theta) = \ones_E(\theta) = \BCAS 0 & \theta\in E \\ \infty & \textnormal{otherwise}. \ECAS
$$
Such subspace constraints arise in many problems, either naturally (\eg{} the constrained lasso by \cite{james2012constrained}) or after reformulation (\eg{} group lasso with overlapping groups).

Before we state our theoretical results, we note that regularizers of the form $\rho(D\theta)$ for some $D\in\reals^{m\times p}$ are geometrically decomposable, as long as $\rho$ is geometrically decomposable. By the geometric decomposability of $\rho,$
\begin{align*}
\rho(D\theta) &= h_A(D\theta) + h_I(D\theta) + h_{E^\perp}(D\theta) \\
&= h_{D^TA}(\theta) + h_{D^TI}(\theta) + h_{D^TE^\perp}(\theta).
\end{align*}
In signal processing, regularizing with $\rho(D\theta)$ for some dictionary $D$ is called \emph{analysis regularization}. We give some examples of M-estimators with geometrically decomposable penalties in Section \ref{sec:main-result}.

\section{Main results}
\label{sec:main-result}

\subsection{Problem setup}
\label{sec:problem-setup}

We begin with a description of the problem at hand. Let $X^{(n)} = \{X_1,\dots,X_n\}$ be $n$ identically distributed observations of some random variable with marginal distribution $\bP.$ We seek to estimate some (unknown) parameters $\theta^\star\in M\subset\reals^p$ of $\bP,$ where $M$ is the \emph{model subspace.} The model subspace is usually low-dimensional and captures the simple structure of the model. For example, $M$ may be the subspace of vectors with a particular support or a subspace of low-rank matrices. We focus on the high-dimensional setting, \ie{} when $n > p.$

Let $\ell$ be a convex and twice-continuously differentiable loss that assigns a cost $\ell(\theta)$ to any parameter $\theta\in E.$ To estimate $\theta^\star$ from the data $X^{(n)},$ we solve the convex optimization problem
\BEQ
\minimize_{\theta\in \reals^p}\,\ell(\theta) + \lambda (h_A(\theta) + h_I(\theta) + h_{E^\perp}(\theta)),
\label{eq:regularized-m-estimator-1}
\EEQ
where $\rho$ is geometrically decomposable in terms of convex sets $A,I\subset\reals^p$ and a subspace $E\subset\reals^p.$ The sets $A,I,E$ are chosen such that $M = E\,\cap\,\linspan(I)^\perp.$ Intuitively, $\linspan(I)\subset M^\perp$ contains deviations from $M$ that we wish to kill. Many common regularized M-estimators possess the decomposable structure given by \eqref{eq:regularized-m-estimator-1}. To gain some intuition, we give three examples, beginning with sparse regression.

\subsubsection{Sparse linear regression}
\label{sec:sparse-regression}

Consider the linear model
\BEQ
y = X\theta + \epsilon,
\label{eq:linear-model}
\EEQ
where $X\in\reals^{n\times p}$ is the design matrix, and $y\in\reals^n$ are the responses. We assume the coefficients $\theta\in\reals^p$ are sparse, \ie{} most of the coefficients are zero. Let $\cS\subset[\,p\,]$ be the support of $\theta,$ and ${\cS^c}$ be the complementary subset of $[\,p\,].$ The model subspace is $\big\{\theta\in\reals^p\mid \theta_{\cS^c} = 0\big\}.$

The \emph{lasso} by \cite{tibshirani1996regression} (also known as \emph{basis pursuit denoising} by \cite{chen2001atomic}) estimates $\thetas$ by the solution of:
\BEQ
\minimize_{\theta\in\reals^p}\frac{1}{2n}\|y - X\theta\|_2^2 + \lambda\|\theta\|_1.
\label{eq:lasso}
\EEQ
The $\ell_1$ norm is geometrically decomposable: $\norm{\theta}_1 = h_{B_{\infty,\cS}}(\theta) + h_{B_{\infty,{\cS^c}}}(\theta),$ where $h_{B_{\infty,\cS}}$ and $h_{B_{\infty,{\cS^c}}}$ are support functions of the sets
\begin{gather*}
B_{\infty,\cS} = \big\{\theta\in\reals^p\mid\norm{\theta}_\infty \le 1,\theta_{\cS^c} = 0\big\} \\
B_{\infty,{\cS^c}} = \big\{\theta\in\reals^p\mid\norm{\theta}_\infty \le 1,\theta_\cS = 0\big\}.
\end{gather*}
It is straightforward to check $\linspan(B_{\infty,{\cS^c}})^\perp = M.$ Thus the lasso possesses the structure given by \eqref{eq:regularized-m-estimator-1}. There is a well-developed theory of the lasso that says, under suitable assumptions on $X,$ the lasso estimator is (consistent and) model selection consistent. In fact, under a stronger beta-min condition, the lasso is \emph{sign consistent.} As we shall see, the aforementioned struc\-ture is the key to the performance of the lasso.

Given an estimate $\htheta,$ there are various ways to assess its performance. We consider two notions: consistency and model selection consistency. An estimate $\hat{\theta}$ is \emph{consistent} (in the $\ell_2$ norm) if the estimation error in the $\ell_2$ norm decays to zero in probability as sample size grows:
$$
\bigl\|\hat{\theta} - \theta^\star\bigr\|_2 \overset{p}{\to} 0\,\text{as }n\to\infty.
$$
An estimate is \emph{model selection consistent} if $\htheta$ is in the \emph{model subspace}:
\BEQ
\Pr(\hat{\theta}\in M) \to 1\,\text{as }n\to\infty.
\label{eq:model-selection-consistent}
\EEQ

\subsection{The main result}

Before we state our main result, we state our assumptions on the problem.
Our two main assumptions are on the (sample) \emph{Fisher information matrix}: $Q = \nabla^2 \ell(\theta^\star).$ The first is \emph{restricted strong convexity (RSC)} and the second is \emph{irrepresentability}.

\begin{assumption}[Restricted strong convexity (RSC)]
\label{asu:rsc}
Let $C\subset\reals^p$ be some (a priori) known convex set containing $\theta^\star.$ The loss function $\ell$ is RSC (on $C\,\cap\, M$) when
\begin{gather}
\Delta^T \nabla^2 \ell(\theta)\Delta \ge m\norm{\Delta}_2^2,\theta \in C\,\cap\, M,\Delta\in (C\,\cap\, M) - (C\,\cap\, M) \\
\|\nabla^2\ell(\theta) - Q\|_2 \le L\norm{\theta -\theta^\star}_2,\theta \in C.
\label{eq:restricted-strong-smoothness}
\end{gather}
for some $m > 0$ and $L < \infty.$
\end{assumption}

The set $C$ is usually taken to be a compact set (see Section \ref{sec:learning-exponential-families} for an example). In these cases, restricted strong smoothness \eqref{eq:restricted-strong-smoothness} holds by the continuity of $\nabla^2\ell.$ Similar notions of restricted strong convexity/smoothness are common in the literature on high-dimensional statistics. For example, the unified framework by \cite{negahban2012unified} requires a (slightly stronger) notion of restricted strong convexity.

For a concrete example, we consider the sparse linear regression problem described in Section \ref{sec:sparse-regression}. When the rows of $X\in\reals^{n\times p}$ are \iid{} Gaussian random vectors with mean zero and covariance $\Sigma,$ \cite{raskutti2010restricted} showed there are constants $m_1,m_2 > 0$ such that
\[
\frac1n\norm{X\Delta}_2^2 \ge m_1\norm{\Delta}_2^2 - m_2\frac{\log p}{n}\norm{\Delta}_1^2\text{ for any }\Delta\in\reals^p
\]
with probability at least $1 - c_1\exp\left(-c_2n\right).$ Their result implies RSC over $\linspan(B_{\infty,\cS})$ with constants $L = 0$ and $m = \frac{m_1}{2}$ as long as $n > 2\frac{m_2}{m_1}\abs{\cS}\log p.$ Thus sparse regression with random Gaussian designs satisfy RSC, even when there are dependencies among the predictors. The result was extended to subgaussian designs by \cite{rudelson2013reconstruction}, also allowing for dependencies among the predictors.

\begin{assumption}[Irrepresentability]
\label{asu:irrepresentable}
There is $\tau\in (0,1)$ such that
\BEQ\textstyle
\sup_{z\,\in\,\partial h_A(M)} V(P_{M^\perp}(QP_M(P_MQP_M)^\dagger P_M z - z))  < 1-\tau,
\label{eq:irrepresentable}
\EEQ
where $V(z) = \inf_y\,\{\gamma_I(y) + \ones_{E^\perp}(z- y)\}$ and $\partial h_A( M) = \bigcup_{\theta\in M}\partial h_A(\theta).$
\end{assumption}

As we shall see, $V$ is a semi-norm: it measures the size of the component of $z$ in $I.$ In particular, $V(z) < 1$ implies
\[
z = z_I + z_{E^\perp}\text{ for some }z_I\in\relint(I)\text{ and }z_{E^\perp}\in E^\perp.
\]

\begin{lemma}
\label{lem:V-semi-norm}
$V$ is finite and sublinear.
\end{lemma}

To interpret the irrepresentable condition, consider again the sparse regression problem. Since $E = \reals^p$ and $Q = \frac1nX^TX,$ \eqref{eq:irrepresentable} simplifies to
\[
\bigl\|X_{\cS^c}^T\bigl(X_{\cS}^T\bigr)^\dagger\sign(\thetas_{\cS})\bigr\|_\infty \le 1 - \tau.
\]
By the properties of support functions, $\partial h_A(M)\subset B_{\infty}.$ Thus it is sufficient to assume
\BEQ
\bigl\|X_{\cS^c}^T\bigl(X_{\cS}^T\bigr)^\dagger\bigr\|_\infty \le 1-\alpha\text{ for some }\alpha\in(0,1).
\label{eq:sparse-regression-irrepresentability}
\EEQ
The rows of $X_{\cS^c}^T\bigl(X_{\cS}^T\bigr)^\dagger$ are the regression coefficients of $x_j,j\in{\cS^c}$ on $X_{\cS}.$ Thus \eqref{eq:sparse-regression-irrepresentability} says the active predictors (columns of $X_\cS$) are not overly well-aligned with the inactive predictors. Ideally, we would like the inactive predictors to be orthogonal to active predictors: $\bigl\|X_{\cS^c}^T\bigl(X_{\cS}^T\bigr)^\dagger\bigr\|_\infty = 0.$ Unfortunately, orthogonality is impossible in the high-dimensional setting. The irrepresentable condition relaxes orthogonality to ``near orthogonality''.

As we shall see, the main result requires the regularization parameter $\lambda$ to be larger than the ``empirical process'' part of the problem. Known results on the convergence rates of regularized M-estimators usually require $\lambda = \Omega(\rho^*(\nabla\ell(\thetas))).$ However, when $\rho$ is not a norm (\eg{} when there are unregularized parameters), $\rho^*(\nabla\ell(\thetas))$ is usually infinite. To allow for unregularized parameters, we relax the requirement to $\lambda  = \Omega(\varrho^*(\nabla\ell(\thetas)))$ for some norm $\varrho$ such that $\rho(\theta) \le \varrho(\theta)$ for any $\theta\in\reals^p.$

Before we state our main result, we describe some constant that appear in the result. Let $B_2$ be the $2$-norm ball. We use $\kappa_\rho$ (resp. $\kappa_{\varrho},\kappa_{\varrho^*}$) to denote the \emph{compatibility constant} between $\rho$ (resp. $\varrho,\varrho^*$) and the $\ell_2$-norm on $M:$
\[\textstyle
\kappa_\rho = \sup_\theta\left\{\rho(\theta)\mid \theta\in B_2\,\cap\, M\right\}
\]
(resp. $\kappa_{\varrho},\kappa_{\varrho^*}$). Similarly, we use $\kappa_{\IC}$ to denote the compatibility constant between the irrep\-resentable term and $\varrho^*:$
\[\textstyle
\kappa_{\IC} = \sup_{\varrho^*(z)\,\le\,1} V(P_{M^\perp}(QP_M(P_MQ P_M)^\dagger P_Mz - z)).
\]
The constants $\kappa_\rho$ and $\kappa_{\IC}$ are finite because $\rho$ and $\varrho^*$ are finite.

\begin{theorem}
\label{thm:main}
Assume $\ell$ and $\rho$ satisfy RSC (on $C\,\cap\, M$) and irrepresentability (Assumptions \ref{asu:rsc} and \ref{asu:irrepresentable}). For any
\BEQ\textstyle
\frac{4\kappa_{\IC}}{\tau}\varrho^*(\nabla\ell(\theta^\star)) < \lambda < \frac{m^2}{2L}\left(2\kappa_{\rho} + \frac{\kappa_\varrho}{\kappa_{\IC}}\frac{\tau}{2}\right)^{-2}\frac{\tau}{\kappa_{\varrho^*}\kappa_{\IC}},
\label{eq:lambda}
\EEQ
the optimal solution to \eqref{eq:regularized-m-estimator-1} is unique,
\BNUM
\item consistent: $\bigl\|\hat{\theta} - \theta^\star\bigr\|_2 \le \frac{2}{m}\Bigl(\kappa_{\rho}+ \frac{\tau}{4}\frac{\kappa_\varrho}{\kappa_{\IC}}\Bigr)\lambda, $
\item model selection consistent: $\hat{\theta} \in M.$
\ENUM
\end{theorem}

Theorem \ref{thm:main} makes a \emph{deterministic} statement about the optimal solution to \eqref{eq:regularized-m-estimator-1}. To use this result to derive consistency and model selection consistency results for a statistical model, we must first verify the loss and penalty satisfies restricted strong convexity/smoothness and irrepresentability. Then, we must select a penalty parameter that satisfies \eqref{eq:lambda} (for some error norm). We know $\varrho^*(\nabla\ell(\theta^\star)) = O_p (\frac{1}{\sqrt{n}})$ for most problems of interest, so, for $n$ large enough, there exist $\lambda$ that satisfies \eqref{eq:lambda}.

\begin{proof}
The proof of Theorem \ref{thm:main} consists of three main steps:
\BNUM
\item Show the solution to a restricted problem \eqref{eq:restricted-problem} is unique and con\-sistent (Lemma \ref{lem:consistency}).
\item Establish a \emph{primal-dual witness (PDW) condition} that ensures all solutions to the original problem are also solutions to the restricted problem (Lemma \ref{lem:dual-certificate}).
\item Construct a primal-dual pair for the original problem from the restricted primal-dual pair that satisfies the dual certificate condition.
\ENUM

Let $(\bar{\theta},\bar{z}_A,\bar{z}_{M^\perp})$ be a primal-dual pair to the restricted problem:
\begin{align}
\minimize_{\theta\in\reals^p}\,\ell(\theta) + \lambda(h_A(\theta) + h_{M^\perp}(\theta)).
\label{eq:restricted-problem}
\end{align}
The restricted primal-dual pair satisfies the first order optimality condition
\begin{gather}
\nabla\ell(\bar{\theta}) +\lambda \bar{z}_A +\lambda \bar{z}_{M^\perp}=0 \\
\bar{z}_A \in \partial h_A (\bar{\theta}),\quad\bar{z}_{M^\perp} \in M^\perp.
\label{eq:kkt-restricted}
\end{gather}
First, we show the solution to the restricted problem is consistent.

\begin{lemma}
\label{lem:consistency}
Assume $\ell$ and $\rho$ satisfy RSC (on $C\,\cap\, M$). For any $\lambda > \frac{4\kappa_{\IC}}{\tau}\varrho^*(P_M\nabla\ell(\theta^\star)),$ the optimal  solution to the restricted problem \eqref{eq:restricted-problem} is uni\-que and consistent: $\bigl\|\bar{\theta} - \theta^\star\bigr\|_2 \le \frac{2}{m}\Bigl(\kappa_{\rho}+ \frac{\tau}{4}\frac{\kappa_\varrho}{\kappa_{\IC}}\Bigr)\lambda.$
\end{lemma}

Next, we establish the PDW condition that ensures all solutions to the original problem are also solutions to the restricted problem.

\begin{lemma}
\label{lem:dual-certificate}
Suppose $\hat{\theta}$ is a primal solution to \eqref{eq:regularized-m-estimator-1}, and  $\hat{z}_A,\hat{z}_I,\hat{z}_{E^\perp}$ are dual solutions; \ie{} $(\hat{\theta},\hat{z}_A,\hat{z}_I,\hat{z}_{E^\perp})$ satisfy
\begin{gather*}
\nabla\ell(\hat{\theta}) +\lambda(\hat{z}_A + \hat{z}_I +\hat{z}_{E^\perp})=0 \\
\hat{z}_I \in \partial h_I (\hat{\theta}),\quad\hat{z}_A \in \partial h_A (\hat{\theta}),\quad\hat{z}_{E^\perp} \in E^\perp.
\end{gather*}
If $\hat{z}_I\in\relint(I)$, then all primal solutions to \eqref{eq:regularized-m-estimator-1} satisfy $h_I(\theta) = 0$.
\end{lemma}

Finally, we use the restricted primal-dual pair to construct a feasible pri\-mal-dual pair for the original problem \eqref{eq:regularized-m-estimator-1}. The optimality conditions of the original problem are
\begin{gather}
\nabla\ell(\hat{\theta}) +\lambda(\hat{z}_A + \hat{z}_I +\hat{z}_{E^\perp})=0 \\
\hat{z}_I \in \partial h_I (\hat{\theta}),\quad\hat{z}_A \in \partial h_A (\hat{\theta}),\quad\hat{z}_{E^\perp} \in E^\perp.
\label{eq:kkt-original}
\end{gather}
Let
\begin{gather*}\textstyle
\hat{z}_I =\argmin_{z} \gamma_I (z)+\ones_{E^\perp}  ( \bar{z}_{M^\perp}-z) \\
\hat{z}_{E^\perp} = \bar{z}_{M^\perp}-\hat{z}_I.
\end{gather*}
The pair $(\bar{\theta},\bar{z}_A,\hat{z}_I,\hat{z}_{E^\perp})$ satisfies \eqref{eq:kkt-original} by construction. Thus $\bar{\theta}$ is a solu\-tion to the original problem. To show $\bar{\theta}$ is the \emph{unique} solution to the original problem, we show $\hat{z}_I$ is PDW feasible: $\hat{z}_I\in\relint(I)$.

The restricted primal-dual pair $(\bar{\theta},\bar{z}_A,\bar{z}_{M^\perp})$ satisfies \eqref{eq:kkt-restricted} and thus the zero reduced gradient condition:
$$
P_M\nabla\ell (\bar{\theta}) +\lambda P_M\bar{z}_A = 0.
$$
We Taylor expand $\nabla\ell$  around $\theta^\star$ (component-wise) to obtain
$$
P_MW + P_MQP_M(\bar{\theta} - \theta^\star) + P_MR + \lambda P_M\bar{z}_A = 0,
$$
where $W = \nabla\ell(\thetas)$ and
$$
R= \nabla \ell(\bar{\theta})-\nabla \ell (\theta^\star) -Q (\bar{\theta}-\theta^\star).
$$
Since $P_MQ P_M$ is invertible on $M$, we solve for $\bar{\theta}$ to obtain
\BEQ
\bar{\theta} = \theta^\star - (P_MQ P_M)^\dagger P_M(W + \lambda\bar{z}_A + R).
\label{eq:solve-for-theta-hat-M}
\EEQ
We Taylor expand $\nabla\ell$ in \eqref{eq:kkt-restricted} around $\theta^\star$ to obtain
$$
W + Q(\bar{\theta} - \theta^\star) + R + \lambda(\bar{z}_A + \bar{z}_{M^\perp}) = 0.
$$
We substitute \eqref{eq:solve-for-theta-hat-M} into this expression to obtain
\BEQ
0 = W - Q (P_MQ P_M)^\dagger P_M(W + \lambda\bar{z}_A + R) + R + \lambda(\bar{z}_A + \bar{z}_{M^\perp}).
\label{eq:substitute-into-solve-for-theta-hat-M-1}
\EEQ
Rearranging, we obtain
\begin{align*}
\bar{z}_{M^\perp} &= \frac{1}{\lambda}\left(Q (P_MQ P_M)^\dagger P_M(W + \lambda \bar{z}_A + R)- W - R - \lambda \bar{z}_A)\right) \\
&= QP_M(P_MQ P_M)^\dagger P_M\bar{z}_A - \bar{z}_A \\
&\pc + \frac{1}{\lambda}\left(QP_M(P_MQ P_M)^\dagger P_M(W + R) - W + R\right).
\end{align*}
Finally, we take $V$'s to obtain
\begin{align}
&V(\bar{z}_{M^\perp}) \le V(P_{M^\perp}(QP_M(P_MQ P_M)^\dagger P_M\bar{z}_A - \bar{z}_A)) \\
&\pc + \frac{1}{\lambda}V(P_{M^\perp}(QP_M(P_MQ P_M)^\dagger W - W))  \\
&\pc+\frac{1}{\lambda}V(P_{M^\perp}(QP_M(P_MQ P_M)^\dagger P_MR - R)).
\label{eq:bound-on-V}
\end{align}
The irrepresentable condition \eqref{eq:irrepresentable} implies the first term is small:
$$
V(P_{M^\perp}(QP_M(P_MQP_M)^\dagger P_M\bar{z}_A - \bar{z}_A)) \le 1-\tau.
$$
Since $V$ is a semi-norm on $M^\perp$, there is some $\kappa_{\IC}$ such that
$$
V(P_{M^\perp}(QP_M(P_MQ P_M)^\dagger W - W))\le\kappa_{\IC}\varrho^*(W).
$$
We substitute these expressions into \eqref{eq:bound-on-V} to obtain
\begin{align*}
V(\bar{z}_{M^\perp}) &\le 1-\tau + \kappa_{\IC}\left(\frac{\varrho^*(W)}{\lambda} + \frac{\varrho^*(R)}{\lambda}\right).
\end{align*}
If we have $\lambda > \frac{4\kappa_{\IC}}{\tau}\varrho^*(W)$, then $\frac{\kappa_{\IC}}{\lambda}\varrho^*(W) \le \frac{\tau}{4}$ and
\begin{align}
V(\bar{z}_{M^\perp}) &< 1-\tau + \frac{\tau}{4} + \frac{\kappa_{\IC}}{\lambda}\rho^*(R).
\label{eq:take-V-2}
\end{align}

\begin{lemma}
\label{lem:R-small}
Assume $\ell$ and $\rho$ satisfy RSC (over $C\,\cap\, M$). For any $\lambda < \frac{m^2}{2L}\left(2\kappa_{\rho} + \frac{\kappa_\varrho}{\kappa_{\IC}}\frac{\tau}{2}\right)^{-2}\frac{\tau}{\kappa_{\varrho^*}\kappa_{\IC}}$, $\frac{\kappa_{\IC}}{\lambda}\varrho^*(R) < \frac{\tau}{4}$.
\end{lemma}

We substitute this bound into \eqref{eq:take-V-2} to obtain
$$
V(\bar{z}_{M^\perp}) < 1-\tau + \frac{\tau}{4} + \frac{\tau}{4} \le 1 - \frac{\tau}{2} <1.
$$
Thus $\hat{z}_I$ is PDW feasible. By Lemma \ref{lem:dual-certificate} and the uniquenss of the solution to the restricted problem, $\btheta$ is also the unique solution to the original problem.
\end{proof}

\subsection{(Partial) converse results}
\label{sec:converse-results}

Although the irrepresentable condition \eqref{eq:irrepresentable} seems cryptic and hard to verify, \cite{zhao2006model} and \cite{wainwright2009sharp} showed it is necessary for sign consistency of the lasso.\footnote{\cite{zhao2006model} and \cite{wainwright2009sharp} refer to the (slightly) stronger condition \eqref{eq:sparse-regression-irrepresentability} as irrepresentability. Thus their converse results are often summarized as irrepresen\-tability is ``almost'' necessary for model selection consistency of the lasso.} In this section, we give necessary conditions for an M-estimator with a geometrically decomposable penalty to be both consistent and model selection consistent.

\begin{theorem}
\label{thm:converse}
Assume $\ell$ and $\rho$ satisfy RSC (on $C\,\cap\, M$) and irrepresentability (Assumptions \ref{asu:rsc} and \ref{asu:irrepresentable}). Further, assume the optimal solution to \eqref{eq:regularized-m-estimator-1} is unique, consistent, and model selection consistent, \ie{}
\[
\htheta\in (\thetas + rB_2)\,\cap\, M.
\]
We must have
\[
\begin{gathered}
\begin{aligned}
&P_{\Mp}QP_M (P_MQ P_M)^\dagger (W + \lambda \hz_A + R) \\
&\pc\in P_{\Mp}(W + R + \lambda (\hz_A + I + E^\perp))  \\
\end{aligned} \\
\hz_A \in\partial h_A((\thetas + rB_2)\,\cap\, M),
\end{gathered}
\]
where $W = \nabla\ell(\theta^\star)$ and $R = \nabla\ell(\hat{\theta}) - W - Q(\hat{\theta} - \theta^\star).$
\end{theorem}

\begin{proof}[Proof of Theorem \ref{thm:converse}]
The proof proceeds like the proof of Theorem \ref{thm:main}. The optimal solution to \eqref{eq:regularized-m-estimator-1} satisfies
\begin{gather}
\nabla\ell(\hat{\theta}) +\lambda(\hz_A + \hz_I +\hz_{E^\perp})=0 \\
\hz_I \in \partial h_I (\hat{\theta}),\quad\hz_A \in \partial h_A (\hat{\theta}),\quad\hz_{E^\perp} \in E^\perp,
\label{eq:first-order-cond}
\end{gather}
Since $\htheta$ is consistent and model selection consistent (by assumption), $\htheta \in (\thetas + rB_2)\,\cap\, M.$ We solve for the error like we did to prove Theorem \ref{thm:main}:
\[
\hat{\theta} - \theta^\star  = - (P_MQ P_M)^\dagger P_M(W + \lambda \hz_A + R).
\]
We plug in the expression for the error to \eqref{eq:first-order-cond} to obtain
\begin{align*}
0 &= W - Q (P_MQ P_M)^\dagger P_M(W + \lambda \hz_A + R) + R + \lambda(\hz_A + \hz_I + \hz_{E^\perp}).
\end{align*}
We project onto $M^\perp$ to obtain the desired result.
\end{proof}

Theorem \ref{thm:converse} is a deterministic statement concerning the solution to \eqref{eq:regularized-m-estimator-1}. It says the (random) term
\BEQ
P_{\Mp}(W + R) - P_{\Mp}QP_M(P_MQP_M)^\dagger(W + R)
\label{eq:random-term}
\EEQ
falls in the set
\BEQ
\begin{aligned}
&P_{\Mp}(\partial h_A((\thetas + rB_2)\,\cap\, M) + I + E^\perp) \\
&\pc- P_{\Mp}QP_M(P_MQP_M)^\dagger \partial h_A((\thetas + rB_2)\,\cap\, M).
\end{aligned}
\label{eq:the-set}
\EEQ
To deduce the necessity of irrepresentability, we must show when irrepresentability is violated, the claims of Theorem \ref{thm:converse} are invalid with positive probability. Although the distribution of \eqref{eq:random-term} is generally hard to characterize, we do not need to completely characterize its distribution. As we shall see, showing the it is symmetric, \ie
$$
\Pr(\eqref{eq:random-term}\in B) = \Pr(\eqref{eq:random-term}\in -B)\text{ for any measurable set }B,
$$
is enough to deduce the necessity of irrepresentability.

%

\begin{corollary}
\label{cor:polyhedral-irrepresentable-necessary}
Assume $\ell$ and $\rho$ satisfy RSC (over $C\,\cap\, M$) and $A$ is a polyhedral set. Further, assume the distribution of \eqref{eq:random-term} is symmetric, and
\[\textstyle
\theta^\star \in \bigcup_{\theta\in\ext(A)}\relint(N_A(\theta)).
\]
When irrepresentability is violated---say
\[\textstyle
\inf_{z\,\in\,\partial h_A(\thetas)} V(P_{M^\perp}(QP_M(P_MQP_M)^\dagger P_M z - z)) \ge 1,
\]
\[\textstyle
\Pr(\htheta\in(\thetas + rB_2)\,\cap\, M) \le \frac12
\]
for any $r$ small enough such that $\thetas + rB_2\subset\bigcup_{x\in\ext(A)}\relint(N_A(x)).$
\end{corollary}

\begin{proof}
Since $\theta^\star \in \bigcup_{x\in\ext(A)}\relint(N_A(x)),$ $\partial h_A(\theta^\star)$ is a point. For any $r$ small enough such that
\[\textstyle
\thetas + rB_2\subset\bigcup_{x\in\ext(A)}\relint(N_A(x)),
\]
$\partial h_A((\thetas + rB_2)\,\cap\, M)$ is also the point $\partial h_A(\theta^\star).$ Thus \eqref{eq:the-set} is given by
\BEQ
P_{\Mp}(\partial h_A(\theta^\star) + I + E^\perp) - P_{\Mp}QP_M(P_MQP_M)^\dagger \partial h_A(\theta^\star).
\label{eq:the-set-2}
\EEQ
When irrepresentability is violated, \eqref{eq:the-set-2} is a convex set that does not contain a relative neighborhood of the origin. Thus there is a halfspace (through the origin) that contains \eqref{eq:the-set-2}. Since the distribution of \eqref{eq:random-term} is symmetric, $\Pr(\eqref{eq:random-term}\in\eqref{eq:the-set}) \le \frac12.$
\end{proof}


\section{Examples}
\label{sec:examples}

We use Theorem \ref{thm:main} to establish the consistency and model selection consistency of the lasso, the generalized lasso, and the regularized maximum likelihood estimator for exponential families in the high-dimensional setting. Our results are nonasymptotic, \ie{} we obtain bounds in terms $n$ and $p$ that hold with high probability.

\subsection{Sparse linear regression}

We return to the sparse linear regression setup described in Section \ref{sec:sparse-regression}. The Fisher information is $\hSigma = \frac1nX^TX.$ We assume
\BNUM
\item RSC (over $\linspan(B_{\infty,\cS})$ and let $C=\reals^p$) and \eqref{eq:sparse-regression-irrepresentability}. Since $\ell$ is quadratic, it satisfies the smoothness condition with $L = 0.$
\item the components of $\epsilon$ are \iid{} subgaussian random variables with mean zero and subgaussian norm $\sigma.$
\ENUM
The assumption \eqref{eq:sparse-regression-irrepresentability} is a stronger condition than irrepresentability. It implies irrepresentability with $\tau = \alpha:$
\BEQ
\begin{aligned}
\bigl\|X_{\cS^c}^T\bigl(X_{\cS}^T\bigr)^\dagger\sign(\thetas_{\cS})\bigr\|_\infty \le \bigl\|X_{\cS^c}^T\bigl(X_{\cS}^T\bigr)^\dagger\bigr\|_\infty\big\|\sign(\thetas)_{\cS}\bigr\|_\infty \le 1-\alpha.
\end{aligned}
\label{eq:sparse-regression-irrepresentability-implies-irrepresentability}
\EEQ

\begin{corollary}
\label{cor:sparse-regression-model-selection-consistency}
Assume $\hSigma$ is RSC (on $\linspan(B_{\infty,\cS})$) and \eqref{eq:sparse-regression-irrepresentability}. For $\lambda$ $ = \frac{8(2 - \alpha)}{\alpha}\sigma\sqrt{\frac{\log p}{n}},$ the lasso estimator is unique,
\BNUM
\item consistent: $\bigl\|\hat{\theta} - \theta^\star\bigr\|_2 \le \frac{4}{m}\Bigl(1+\frac{4(2-\alpha)}{\alpha}\Bigr)\sigma\sqrt{\frac{\abs{\cS}\log p}{n}}, $
\item model selection consistent: $\hat{\theta}_{\cS^c} = 0$
\ENUM
with probability at least $1 - 2p^{-1}.$

Further, if $\min_{a\in\cS}\abs{\theta^\star_a} > \frac{4}{m}\Bigl(1+\frac{4(2-\alpha)}{\alpha}\Bigr)\sigma\sqrt{\frac{\abs{\cS}\log p}{n}},$ then the lasso estimator is also sign consistent: $\sign(\hat \theta_\cS) = \sign(\thetas_\cS).$
\end{corollary}

\begin{proof}
Before we apply Theorem \ref{thm:main}, we compute the constants $\kappa_\rho,\kappa_\varrho$ and $\kappa_{\IC}.$
Since the regularizer is finite (it's a norm), its dual semi-norm is finite. To keep things simple, we let $\varrho = \norm{\cdot}_1.$ The constant $\kappa_\rho = \kappa_\varrho$ is
\[\textstyle
\kappa_\rho = \sup_\theta\left\{\norm{\theta}_1\mid \theta\in B_2\,\cap\,\linspan(B_{\infty,\cS})\right\} = \sqrt{\abs{\cS}}.
\]
Similarly, the constant $\kappa_{\IC}$ is given by
\begin{align*}
&\bigl\|P_{B_{\infty,{\cS^c}}}(\hSigma P_{B_{\infty,\cS}}(P_{B_{\infty,\cS}}\hSigma P_{B_{\infty,\cS}})^\dagger P_{B_{\infty,\cS}}z - z)\bigr\|_\infty \\
&\pc\le\bigl\|X_{\cS^c}^T\bigl(X_{\cS}^T\bigr)^\dagger z_{\cS}\bigr\|_\infty + \norm{z_{\cS^c}}_\infty \le (2-\alpha)\norm{z}_\infty.
\end{align*}
is at most $2 - \alpha.$

To apply Theorem \ref{thm:main}, we check $\lambda = \frac{8(2-\alpha)}{\alpha}\sigma\sqrt{\frac{\log p}{n}}$ satisfies the assumptions. Since the loss function is quadratic, it satisfies the smoothness condition \eqref{eq:restricted-strong-smoothness} with $L = 0.$ Thus any $\lambda < \infty$ satisfies the upper bound in \eqref{eq:lambda}. We check our choice also satisfies the lower bound in \eqref{eq:lambda}. By \cite{vershynin2010introduction}, Proposition 5.10 and a union bound,
$$\textstyle
\Pr\left(\big\|\nabla\ell(\theta^\star)\big\|_{\infty} > t \right) = \Pr\big(\norm{X^T \epsilon}_\infty > nt \big) \le 2\exp\Bigl(-\frac{nt^2}{2 \sigma^2} +\log p\Bigr).
$$
Thus
\begin{align*}
&\Pr\biggl( \frac{4(2-\alpha)}{\alpha}\big\|\nabla\ell(\theta^\star)\big\|_\infty > \frac{8(2-\alpha)}{\alpha}\sigma \sqrt{\frac{\log p}{n}}\,\biggr) \\
&\pc\le 2\exp(-2\log p + \log p) = 2 p^{-1}.
\end{align*}
Consequently, the claims of Theorem \ref{thm:main} are valid with probability at least $1 - 2 p^{-1}:$
\BNUM
\item $\bigl\|\hat{\theta} - \theta^\star\bigr\|_2 \le \frac{2}{m}\Bigl(1 + \frac{\alpha}{4(2 - \alpha)}\Bigr)\sqrt{\abs{\cS}}\lambda = \frac{4}{m}\Bigl(1+\frac{4(2-\alpha)}{\alpha}\Bigr)\sigma\sqrt{\frac{\abs{\cS}\log p}{n}},$
\item $\htheta \in \linspan(B_{\infty,\cS}) = \{\theta\in\reals^p\mid \theta_{\cS^c} = 0\}.$
\ENUM

An easy consequence of (1) is $\bigl\|\htheta_a - \thetas_a\bigr\|_\infty \le \theta_{\min}.$ Thus $\htheta$ is sign consis-tent: $\sign(\hat \theta_a) = \sign(\theta^\star _a)$ for any $a\in\cS$ such that $\abs{\theta^\star_a} > \theta_{\min}.$
\end{proof}

As we saw in Section \ref{sec:geometrically-decomposable}, analysis regularizers of the form $\rho(D\theta)$ are geometrically decomposable. A prominent example of $\ell_1$ analysis regularization is the \emph{generalized lasso:}
\BEQ
\minimize_{\theta\in\reals^p}\frac{1}{2n}\norm{y-X\theta}_2^2 + \lambda\norm{D\theta}_1.
\label{eq:generalized-lasso}
\EEQ
The underlying (statistical) model is a straightforward modification of the linear model \eqref{eq:linear-model}: we assume $D\thetas$ (instead of $\thetas$) is sparse. The sparsity of $D\theta$ usually translates to some desirable structural or geometric property of $\theta.$ We refer to Section 2 in \cite{tibshirani2011solution} for some examples.

The model subspace is $\{\theta\in\reals^p\mid (D\theta)_{\cS^c} = 0\},$ where $\cS\subset[\,m\,]$ is the support of $D\theta.$ It's straightforward to check $\linspan(D^TB_{\infty,\cS^c})^\perp = M.$ Thus the generalized lasso possesses the structure given by \eqref{eq:regularized-m-estimator-1}. To study the model selection properties of the generalized lasso, we assume
\BEQ
\bigl\|D_{\cS^c}X^T\bigl(D_{\cS}X^T\bigr)^\dagger\sign(\thetas_{\cS})\bigr\|_\infty \le 1 - \tau
\label{eq:ic}
\EEQ
in lieu of \eqref{eq:sparse-regression-irrepresentability}. The assumption \eqref{eq:ic} is equivalent to irrepresentability. It is usually referred to as an \emph{identifiability criterion (IC)}. Given IC \eqref{eq:ic}, we derive the analog of Corollary \eqref{cor:sparse-regression-model-selection-consistency} for the generalized lasso.

\begin{corollary}
Assume $\hSigma$ is RSC (on $\linspan(D^TB_{\infty,\cS^c})^\perp$) and \eqref{eq:ic}. For $\lambda$ $ = \frac{8\kappa_{\IC}}{\tau}\sigma\sqrt{\frac{\log p}{n}},$ the generalized lasso estimator is unique,
\BNUM
\item consistent: $\bigl\|\hat{\theta} - \theta^\star\bigr\|_2 \le \frac{4}{m}\Bigl(\kappa_{\varrho} + \frac{4\kappa_{\IC}}{\tau}\kappa_{\rho}\Bigr)\sigma\sqrt{\frac{\log p}{n}},$
\item model selection consistent: $(D\theta)_{\cS^c} = 0$
\ENUM
with probability at least $1 - 2p^{-1}.$
\end{corollary}

\begin{proof}
Before we apply Theorem \ref{thm:main}, we compute the constants $\kappa_\rho$ and $\kappa_\varrho.$ When $D$ has a (nontrivial) null space, the regularizer is not a norm. To allow for the possibility, we let $\varrho = \norm{\cdot}_1$ and set $\lambda > \frac{4\kappa_{\IC}}{\tau}\|\nabla\ell(\thetas)\|_\infty.$ The constants $\kappa_{\rho},\kappa_{\varrho}$ are
\begin{gather*}
\textstyle\kappa_\rho = \sup_\theta\left\{\norm{D\theta}_1\mid \theta\in B_2\,\cap\,\linspan(D^TB_{\infty,\cS^c})^\perp\right\} \\
\textstyle\kappa_{\varrho} = \sup_\theta\left\{\norm{\theta}_1\mid \theta\in B_2\,\cap\,\linspan(D^TB_{\infty,\cS^c})^\perp\right\}.
\end{gather*}

By an argument similar to the argument in the proof of Corollary \ref{cor:sparse-regression-model-selection-consistency}, $\lambda = \frac{8\kappa_{\IC}}{\tau}\sigma\sqrt{\frac{\log p}{n}}$ satisfies the assumptions of Theorem \ref{thm:main} with probability at least $1 - 2 p^{-1}.$ Thus, with probability at least $1 - 2 p^{-1},$
\BNUM
\item $\bigl\|\hat{\theta} - \theta^\star\bigr\|_2 \le \frac{2}{m}\Bigl(\kappa_{\rho}+ \frac{\tau}{4}\frac{\kappa_\varrho}{\kappa_{\IC}}\Bigr)\lambda = \frac{2}{m}\Bigl(2\kappa_{\varrho} + \frac{8\kappa_{\IC}}{\tau}\kappa_{\rho}\Bigr)\sigma\sqrt{\frac{\log p}{n}},$
\item $\htheta \in \linspan(D^TB_{\infty,\cS^c})^\perp = \{\theta\in\reals^p\mid (D\theta)_{\cS^c} = 0\}.$
\ENUM
\end{proof}

\subsection{Learning exponential families}
\label{sec:learning-exponential-families}

We turn our attention to a problem with a non-quadratic loss function. Recall an exponential family is a distribution of the form
$$
\Pr(x;\theta) = h(x) \exp\left(\theta^T\phi(x) - A(\theta)\right),
$$
where $\theta$ are the \emph{natural parameters}, $\phi(X)\in\reals^p$ are \emph{sufficient statistics}. We assume $\thetas$ is \emph{group-sparse}, \ie{} the components of $\thetas$ are organized in (possibly overlapping) groups and only a few groups are active. Let $\cG$ be the collection of groups and $\cS$ be the subset of active groups. The model subspace is $M = \{\theta\in\reals^p\mid\theta_g = 0\text{ for any }g\in\cS^c\}.$

Given independent observations $X^{(n)} = \{X_1,\dots,X_n\},$ we seek to estimate $\thetas$ by the regularized \emph{maximum likelihood estimator (MLE):}
\BEQ
\minimize_{\theta\,\in\,E\,\subset\,\reals^p}-\frac1n\sum_{i=1}^n\phi(x^{(i)})^T\theta + A(\theta) + \lambda\norm{\theta}_{2/1}.
\label{eq:regularized-mle}
\EEQ
The $\ell_1/\ell_2$ norm is geometrically decomposable:
$$
\norm{\theta}_{2/1} = \sum_{g\in\cG}\norm{\theta_g}_2 = h_{B_{\infty/2,\cS}}(\theta) + h_{B_{\infty/2,\cS^c}}(\theta),
$$
where $h_{B_{\infty/2,\cS}}$ and $h_{B_{\infty/2,\cS^c}}$ are support functions of the sets
\begin{gather*}
\textstyle B_{\infty/2,\cS} = \left\{\theta\in\reals^p\mid\max_{g\in\cG}\norm{\theta_g}_2 \le 1,\theta_g = 0\text{ for any }g\in {\cS^c}\right\} \\
\textstyle B_{\infty/2,\cS^c} = \left\{\theta\in\reals^p\mid\max_{g\in\cG}\norm{\theta_g}_2 \le 1,\theta_g = 0\text{ for any }g\in \cS\right\}.
\end{gather*}
It is easy to check $\linspan(B_{\infty/2,\cS}) = M.$ Thus \eqref{eq:regularized-mle} has the structure given by \eqref{eq:regularized-m-estimator-1}. The Fisher information is $Q = \nabla^2 A(\thetas).$ We assume $Q$ satisfies RSC (over $\linspan(B_{\infty/2,\cS})$) and irrepresentability.



First, we establish two auxiliary results: (i) a concentration result for $W$ and (ii) the optimal solution to \eqref{eq:regularized-mle} is contained in some compact subset of the model subspace.

\begin{lemma}
\label{lem:subexp-concentration}
The random variable $W$ satisfies
\begin{align*}
&\Pr \left(\left|\frac{\partial \ell}{\partial \theta_j} (\theta^\star)\right| >t \right) \le 2\exp\left( -cn \min\left(\frac{t^2} { \max_{g \in \cG} |g| K ^2}, \frac{t}{ \max_{g \in \cG} \sqrt{|g|}K}\right)\right)\\
&\Pr \left(\max_{g \in \cG} \norm{\nabla_{\theta_g} \ell (\theta^\star) }_2 >t \right)\\ &\le2\exp\left(\log |\cG| -cn \min\left(\frac{t^2} { \max_{g \in \cG} |g| K ^2}, \frac{t}{ \max_{g \in \cG} \sqrt{|g|}K}\right)\right)\end{align*}
for some absolute constant $c>0$ and a constant $K$ that is independent of  $n$.
\end{lemma}

\begin{lemma}
\label{lem:exp-family-compact-subset}
The optimal solution to \eqref{eq:regularized-mle} satisfies
\begin{align*}
 \bigl\|\hat{\theta}\bigr\|_{2/1}  \le\frac{1}{(\lambda-\norm{ \phi^n - \phi^\star}_{2,\infty}) }\left(   \lambda \norm{ \theta^\star}_{2/1}+ \norm{ \phi^n - \phi^\star}_{2,\infty}\norm{\theta^{\star}}_{2/1}\right)\\
 A(\hat \theta) \le \norm{\theta^{\star}}_{2/1} \norm{\phi^n}_{2,\infty} + \bigl\|\hat{\theta}\bigr\|_{2/1} \norm{\phi^n}_{2,\infty} +A(\theta^\star) +\lambda \norm{\theta^\star}_{2/1}
\end{align*}

where $\phi^n = \frac{1}{n} \sum_{i=1}^n \phi(x^{(i)}) $ and $\phi^\star = \Expect_{\theta^\star}[ \phi(X)]$.
\end{lemma}

We use these two results to establish the consistency and model selection consistency of the regularized MLE.

\begin{corollary}
\label{cor:regularized-mle-consistent}
Suppose we are given samples $x^{(1)},\dots,x^{(n)}$ drawn \iid\ from a regular  \emph{exponential family} with unknown parameters $\theta^\star$, $\kappa_{\IC} \ge \tau$, and Assumption \ref{asu:irrepresentable} is satisfied.  Select
$$
\lambda = \frac{3 \kappa_{\IC}}{\tau} \max_{g \in G} \sqrt{|g|} K \sqrt{\frac{\log |\cG|}{cn}}
$$
and the sample size
$$
n> \max \left( \frac{36}{c}\frac{ \kappa_{\IC}^4}{\tau^4}\max_{g \in G} |g| K^2\log |\cG|\frac{L^2}{m^4}\left(2\sqrt{|\cS|} + \frac{\tau}{2\kappa_{\IC}}\sqrt{|\cS|}\right)^4, \left(\frac{3}{2}\right)^2 \frac{\log |\cG|}{c}\right)
$$
where $C:=\left\{\theta\ \big| \norm{\theta}_{2/1} \le 4 \norm{\theta^\star}_{2/1} \text{ and }A(\theta)\le R\right\}$, $c>0$ is an absolute constant, and $K$ is a constant independent of $n$ defined in Lemma \ref{lem:suff-stat-subexp}.
With probability at least $1-2 |\cG|^{-5/4}$, the optimal solution to \eqref{eq:regularized-mle} satisfies
\BNUM
\item $\bigl\|\hat{\theta} - \theta^\star\bigr\|_2 \le \frac{6}{m}   \frac{\kappa_{\IC}}{\tau} \max_{g \in G} \sqrt{|g|} K\left(\sqrt{|\cS|} + \frac{\tau}{2\kappa_{\IC}}\sqrt{|\cS|}\right) \sqrt{\frac{\log |\cG|}{cn}}$
\item $\hat{\theta}_g = 0,g\in{\cS^c}$.
\ENUM
Furthermore if we assume the beta-min condition \[
\norm{\theta^\star_g}_2 > \frac{6}{m}   \frac{\kappa_{\IC}}{\tau} \max_{g \in G} \sqrt{|g|} K\left(\sqrt{|\cS|} + \frac{\tau}{2\kappa_{\IC}}\sqrt{|\cS|}\right) \sqrt{\frac{\log |\cG|}{cn}}\] for all $g \in \cS$, then all groups  $g\in \cS$ are correctly estimated as non-zero, $\|\hat \theta _g\|_2 > 0$.
\label{cor:mle-consistent}
\end{corollary}

\section{Model selection properties of regularized M-estimators with weakly decomposable penalties}
\label{sec:weakly-decomposable}

\subsection{Background and problem setup}
\label{sec:low-rank-multivariate-regression}

Geometric decomposability, although general, excludes some common regularizers. An important example is the nuclear norm:
\[\textstyle
\norm{\Theta}_* = \sum_{j=1}^{r}\sigma_j(\Theta),
\]
where $r$ is the rank of $\Theta\in\reals^{p_1\times p_2}$ and $\sigma_j(\Theta),j=1,\dots,r$ are its singular values. The motivating example we have in mind is low-rank multivariate regression. Consider the (multivariate) generalization of the linear model:
\BEQ
Y = X\Thetas + W,
\label{eq:multivariate-linear-model}
\EEQ
where the rows of $Y\in\reals^{n\times p_2}$ are (multivariate) responses. We assume the matrix of coefficients $\Thetas\in\reals^{p_1\times p_2}$ has rank $r \ll\min\{p_1,p_2\}.$
Given observations $\left\{(x_i,y_i)\right\}_{i=1}^n,$ a standard approach to estimating the unknow $\Thetas$ is \emph{nuclear norm minimization:}
\BEQ
\minimize_{\Theta\in\reals^{p_1\times p_2}}\,\frac{1}{2n}\norm{Y - X\Theta}_{\Fro}^2 + \lambda\norm{\Theta}_*.
\label{eq:low-rank-multivariate-regression}
\EEQ
\cite{bach2008consistency} showed that nuclear norm minimization is \emph{rank consistent}, \ie
\BEQ
\Pr\bigl(\rank(\hTheta) = r\bigr) \to 1\text{ as }n\to\infty,
\label{eq:rank-consistent}
\EEQ
subject to irrepresentability. Although rank consistency does not fit into our notion of model selection consistency because the set of rank $r$ matrices is not a subspace, our results may be used to derive a non-asymptotic form of Bach's rank consistency result.

To study the rank consistency of nuclear norm minimization, we consider an alternative notion of decomposability: \emph{weak decomposability.}

\begin{definition}
\label{def:weakly-decomposable}
A regularizer is \emph{weakly decomposable} at $\thetas\in\reals^p$ in terms of convex sets $A,I\subset\reals^p$ if it is sublinear and
\begin{gather*}
\partial\rho(\thetas) = \partial h_A(\thetas) + \partial h_I(\thetas).
\end{gather*}
We assume $A,I$ are bounded and $0\in\relint(I).$
\end{definition}

Weak decomposability is more general than geometric decomposability. However, the structure of the subdifferential of a weakly decomposable pe\-nalty at $\thetas$ is very similar to that of a geometrically decomposable penalty. Consequently, the directional derivative of $\rho$ at $\thetas$ along $\Delta$ is geometrically decomposable:
\[
\brho(\thetas,\Delta) = h_{\partial h_A(\thetas)}(\Delta) + h_{\partial h_I(\thetas)}(\Delta).
\]
As we shall see, the geometric decomposability of $\brho(\thetas,\Delta)$ is the key to the model selection properties of weakly decomposable penalties.

The problem setup is similar to the setup in Section \ref{sec:problem-setup}. To keep things simple, we focus on regularized least squares. Given $n$ identically distri\-buted observations of some random variable, we estimate some parameters $\thetas\in M\subset\reals^p$ of its distribution by
\BEQ
\minimize_{\theta\in \reals^p}\,\frac12\theta^TQ\theta - q^T\theta + \lambda \rho(\theta),
\label{eq:regularized-m-estimator-2}
\EEQ
where $\rho$ is weakly decomposable in terms of convex sets $A,I\subset\reals^p.$ The sets $A,I$ are chosen such that $M = \linspan(I)^\perp.$ As we shall see, the nuclear norm minimization problem has the form given by \eqref{eq:regularized-m-estimator-2}.

\subsection{Dual consistency of regularized M-estimators}

To study the model selection properties of \eqref{eq:regularized-m-estimator-2}, we compare the its optimal solution to the optimal solution to a linearized problem
\BEQ
\minimize_{\theta\in\reals^p}\,\frac12\theta^TQ\theta - q^T\theta  + \lambda(\rho(\thetas) + \brho(\thetas,\theta - \thetas)).
\label{eq:linearized-problem}
\EEQ
Since the objective functions of \eqref{eq:linearized-problem} and \eqref{eq:regularized-m-estimator-2} are similar, we expect the (optimal) solutions are close. Unfortunately, due to the lack of strong convexity, we cannot conclude the solutions are close. However, as we shall see, the dual solutions are close.

After a change of variables, the linearized problem is
\BEQ
\minimize_{\Delta\in\reals^p}\,\frac12\Delta^TQ\Delta + (Q\thetas - q)^T\Delta + \lambda( h_{\partial h_A(\thetas)}(\Delta) + h_I(\Delta)),
\label{eq:linearized-problem-cov}
\EEQ
We recognize \eqref{eq:linearized-problem-cov} possesses the decomposable structure given by \eqref{eq:regularized-m-estimator-1}. By Theorem \ref{thm:main}, a primal-dual pair $(\bDelta,\bz_A,\bz_I)$ that satisfies
\BEQ
\begin{gathered}
Q(\thetas + \bDelta) - q + \lambda(\bz_A + \bz_I) = 0 \\
\bz_A\in\partial h_A(\thetas),\quad\bz_I \in I
\end{gathered}
\label{eq:linearized-problem-cov-kkt}
\EEQ
is unique. Further, $\bDelta$ is consistent, and $\bz_I$ is PDW feasible. We summarize the properties of $(\bDelta,\bz_A,\bz_I)$ in a lemma.

\begin{lemma}
\label{lem:linearized-problem-consistency}
Assume $Q$ and $\brho$ satisfy RSC (on $\linspan(I)^\perp$) and irrepresentability. For any $\lambda > \frac{4\kappa_{\IC}}{\tau}\bvarrho^*(Q\thetas -q),$ the unique primal-dual pair for \eqref{eq:linearized-problem-cov} $(\bDelta,\bz_A,\bz_I)$ is
\BNUM
\item consistent: $\bigl\|\bDelta\bigr\|_2 \le \frac{2}{m}\Bigl(\kappa_{\brho}+ \frac{\tau}{4}\frac{\kappa_{\bvarrho}}{\kappa_{\IC}}\Bigr)\lambda.$
\item PDW feasible: $\gamma_I(\bz_I) \le 1 - \frac{\tau}{2}.$
\ENUM
\end{lemma}

The main result shows the dual solutions to \eqref{eq:regularized-m-estimator-2} and \eqref{eq:linearized-problem-cov} are close.

\begin{theorem}
\label{thm:dual-consistency}
Assume $Q$ and $\brho$ satisfy RSC (on $\linspan(I)^\perp$) and irrepresentability. For any $\lambda > \frac{4\kappa_{\IC}}{\tau}\bvarrho^*(Q\thetas - q),$ the optimal dual solutions to \eqref{eq:linearized-problem} and \eqref{eq:regularized-m-estimator-2} satisfy
\[\textstyle
\norm{\bz_A + \bz_I - \hz}_2^2 \le \frac{\norm{Q}_2}{\lambda}\bigl(R(\bDelta) - R(\hDelta)\bigr),
\]
where $R(\Delta) = \rho(\thetas + \Delta) - \rho(\thetas) - \brho(\thetas,\Delta).$
\end{theorem}

\begin{proof}
After a change of variables, the original problem is
\BEQ
\minimize_{\Delta\in\reals^p}\,\frac12\Delta^TQ\Delta + (Q\thetas - q)^T\Delta + \lambda\rho(\thetas + \Delta).
\label{eq:regularized-m-estimator-2-cov}
\EEQ
Its optimality conditions are
\BEQ
\begin{gathered}
Q(\thetas + \hDelta) - \gamma + \lambda\hz = 0 \\
\hz\in \partial\rho(\thetas + \hDelta).
\end{gathered}
\label{eq:regularized-m-estimator-2-cov-kkt}
\EEQ
Let $\bDelta$ and $\hDelta$ be the solutions to \eqref{eq:linearized-problem-cov} and \eqref{eq:regularized-m-estimator-2-cov}. By Fermat's rule, $\bz_A + \bz_I$ and $\hz_A + \hz_I$ are also the dual solutions to \eqref{eq:linearized-problem} and \eqref{eq:regularized-m-estimator-2}. We subtract \eqref{eq:regularized-m-estimator-2-cov-kkt} from \eqref{eq:linearized-problem-cov-kkt} to obtain
\BEQ
Q(\hDelta - \bDelta) = \lambda(\bz_A + \bz_I - \hz).
\label{eq:dual-consistency-1}
\EEQ
To complete the proof, we show $\bigl\|Q(\hDelta - \bDelta)\bigr\|_2^2$ is small.
By inspection of the optimality conditions \eqref{eq:regularized-m-estimator-2-cov-kkt} and \eqref{eq:linearized-problem-cov-kkt}, $\bDelta$ and $\hDelta$ are also the solutions to
\[
\begin{gathered}
\minimize_{\Delta\in E} \bDelta^T Q\Delta + (Q\thetas - q)^T\Delta + \lambda\brho(\thetas,\Delta) \\
\minimize_{\Delta\in E} \hDelta^T Q\Delta + (Q\thetas - q)^T\Delta + \lambda(\rho(\thetas + \Delta) - \rho(\thetas)).
\end{gathered}
\]
Since $\bDelta$ and $\hDelta$ are their respective optimal solutions, we know
\begin{gather*}
\begin{aligned}
&\bDelta^T Q\bDelta + (Q\thetas - q)^T\bDelta + \lambda\brho(\thetas,\bDelta) \\
&\pc\le \bDelta^T Q\hDelta + (Q\thetas - q)^T\hDelta + \lambda\brho(\thetas,\hDelta),
\end{aligned}
\\
\begin{aligned}
&\hDelta^T Q\hDelta + (Q\thetas - q)^T\hDelta + \lambda(\rho(\thetas + \hDelta) - \rho(\thetas)) \\
&\pc\le \hDelta^T Q\bDelta + (Q\thetas - q)^T\bDelta + \lambda\left(\rho(\thetas + \bDelta) - \rho(\thetas)\right).
\end{aligned}
\end{gather*}
We add the inequalities and rearrange to obtain
\[
(\bDelta - \hDelta)^TQ(\bDelta - \hDelta) = \|\Delta\|_Q^2 \le \lambda\bigl(R(\bDelta) - R(\hDelta)\bigr).
\]
where $R(\Delta) = \rho(\thetas + \Delta) - \rho(\thetas) - \brho(\thetas,\Delta).$ Since $\|Q\Delta\|_2^2 \le \norm{Q}_2\|\Delta\|_Q^2,$
\[
\bigr\|Q(\hDelta - \bDelta)\bigr\|_2^2 \le \norm{Q}_2\bigr\|\hDelta - \bDelta\bigr\|_Q^2 \le \norm{Q}_2\lambda\bigl(R(\bDelta) - R(\hDelta)\bigr).
\]
We plug in \eqref{eq:dual-consistency-1} to reach the stated conclusion.
\end{proof}

\subsection{Rank consistency of low-rank multivariate regression}

We return to the low-rank multivariate regression problem described in Section \ref{sec:low-rank-multivariate-regression}. The nuclear norm is weakly decomposable. Let $\Thetas = U\Sigma V^T$ be the (full) SVD of $\Thetas$ and define the sets
\begin{align*}
A &= \left\{\Theta\in B_{\op}\subset\reals^{p_1\times p_2}\mid\Theta = U_rDV_r^T\text{ for some diagonal }D\right\} \\
I &= \left\{\Theta\in B_{\op}\subset\reals^{p_1\times p_2}\mid\Theta = U_{p_1-r}DV_{p_2-r}^T\text{ for some diagonal }D\right\},
\end{align*}
where $U_r,U_{p_1-r}$ (resp. $V_r,V_{p_2-r}$) are submatrices of $U$ (resp. $V$) consisting of the top $r$ and bottom $p_1-r$ left (resp. $p_2 - r$ right) singular vectors of $\Thetas.$ It's straightforward to check the nuclear norm is weakly decomposable  at $\Thetas$ in terms of $A,I.$ Since $A + I\subset B_{\op},$
\[
\norm{\Theta}_* = h_{B_{\op}}(\Theta) \ge h_A(\Theta) + h_I(\Theta).
\]

Before we delve into the rank consistency of low-rank multivariate regression, we state our assumptions on the problem. To keep notation manageable, we adopt operator theoretic notation. Let $\vec{X}\in\reals^{p_1p_2}$ be the vectorized form of $X\in\reals^{p_1\times p_2}.$ In operator notation, the model is
\BEQ
\vec{Y} = \cX(\Thetas) + \vec{W},
\label{eq:matrix-linear-model}
\EEQ
where $\cX:\reals^{p_1\times p_2}\to\reals^n$ is a linear operator. Since $\cX$ is linear, we abuse notation by writing $\cX\Theta = \cX(\Theta).$ The Fisher information $Q:\reals^{p_1\times p_2}\to\reals^{p_1\times p_2}$ is given by $\frac1n \cX^*\cX.$ We assume
\BNUM
\item RSC and
\BEQ\textstyle
\sup_{Z\,\in\,B_{\op}}\norm{U_{p_1-r}^T\bigl[P_IQP_{I^\perp}(P_{I^\perp}QP_{I^\perp})^\dagger Z\bigr]V_{p_2-r}}_{\op} \le 1 - \alpha,
\label{eq:low-rank-multivariate-regression-irrepresentability}
\EEQ
where $P_I:\reals^{p_1\times p_2}\to\reals^{p_1\times p_2}$ (resp. $P_{I^\perp}$) is the projector onto $\linspan(I)$ (resp. $\linspan(I)^\perp$).
\item the entries of $W$ are \iid{} subgaussian random variables with mean zero and subgaussian norm $\sigma.$
\ENUM
The assumption \eqref{eq:low-rank-multivariate-regression-irrepresentability} is stronger than irrepresentability. It implies irrepresentability with $\tau = \alpha:$
\BEQ
\begin{aligned}
&\textstyle\norm{U_{p_1-r}^T\bigl[P_I\bigl(QP_{I^\perp}(P_{I^\perp}QP_{I^\perp})^\dagger U_rV_r^T - U_rV_r^T\bigr)\bigr]V_{p_2-r}}_{\op} \\
&\textstyle\pc= \norm{U_{p_1-r}^T\bigl[P_I\bigl(QP_{I^\perp}(P_{I^\perp}QP_{I^\perp})^\dagger U_rV_r^T\bigr)\bigr]V_{p_2-r}}_{\op} \\
&\textstyle\pc\le\sup_{Z\,\in\,B_{\op}}\norm{U_{p_1-r}^T\bigl[P_IQP_{I^\perp}(P_{I^\perp}QP_{I^\perp})^\dagger Z\bigr]V_{p_2-r}}_{\op}.
\end{aligned}
\label{eq:low-rank-multivariate-regression-irrepresentability-implies-irrepresentability}
\EEQ
We make the stronger assumption to obtain an explicit expression for the constant $\kappa_{\IC}$ (in terms of the constant $\alpha$).

The final ingredient we require is a ``Taylor's theorem'' for the nuclear norm that says the nuclear norm is well-approximated by its linearization.

\begin{lemma}
\label{lem:nuclear-norm-taylors-thm}
Let $s_r$ be smallest nonzero singular value of $\Thetas.$ For any $\Delta\in\linspan(I)^\perp,\norm{\Delta}_{\op} < \frac{s_r}{2},$ we have
\[
\norm{\Thetas + \Delta}_* - \norm{\Thetas}_* - \tr\bigl(U_r^T\Delta V_r\bigr) \le \frac{4}{3s_r}\norm{\Delta}_{\Fro}^2.
\]
\end{lemma}
We put the pieces togther to conclude low-rank multivariate regression is rank consistent.

\begin{corollary}
\label{cor:nuclear-norm-minimization-rank-consistency}
Assume $Q$ is RSC (over $\linspan(I)^\perp$) and \eqref{eq:low-rank-multivariate-regression-irrepresentability}. For $\lambda$ $= \frac{8(2-\alpha)}{\alpha}\sigma\sqrt{\frac{p_1+p_2}{n}},$ the optimal solution to \eqref{eq:low-rank-multivariate-regression} is unqiue and rank consistent when
\[
n > \max\left\{\frac{128^2}{9s_r^2}\frac{M^2}{m^4}\frac{(\sqrt{2} + \alpha')^4}{\alpha^4\alpha'^2} r^2,\frac{16}{m^2}(\sqrt{2} + \alpha')^2r\right\}\sigma^2(p_1+p_2)
\]
with probability at least $1 - c_1e^{-c_2(p_1 + p_2)}.$ The constants $M$ and $\alpha'$ are given by $\sup_{\norm{\Delta}_{\Fro}\,\le\,1}\norm{Q\Delta}_{\Fro}$ and $\frac{4(2-\alpha)}{\alpha}.$
\end{corollary}

\begin{proof}
To show $\hTheta$ has rank at most $r,$ it suffices to show the optimal dual solution $\hU\hV^T$ has no more than $r$ unit singular values. At a high level, the proof consists of three steps:
\BNUM
\item Show the unique (feasible) primal-dual pair to a linearized problem $\left(\bDelta, U_rV_r^T, \bU_{p_1 - r}\bV_{p_2 - r}^T\right)$ is consistent and PDW feasible (Lemma \ref{lem:linearized-low-rank-multivariate-regression-consistency}).
\item Invoke Theorem \ref{thm:dual-consistency} to show $\hU\hV^T$ is close to the optimal dual solution to the linearized problem $U_rV_r^T + \bU_{p_1 - r}\bV_{p_2 - r}^T.$ Since $\bU_{p_1 - r}\bV_{p_2 - r}^T$ is PDW feasible, its singular values are bounded away from one.
\item Apply a singular value perturbation result to conclude $\hU\hV^T$ has (no more than) $r$ unit singular values.
\ENUM

Consider the linearized problem
\BEQ
\minimize_{\Delta\in\reals^{p_1\times p_2}}\,\frac{1}{2n}\norm{Y - X(\Thetas + \Delta)}_{\Fro}^2 + \lambda\left(\tr\bigl(U_r^T\Delta V_r\bigr) + \|U_{p_1-r}^T\Delta V_{p_2-r}\|_{*} \right).
\label{eq:nuclear-norm-minimization-linearized-cov}
\EEQ
We apply Lemma \ref{lem:linearized-problem-consistency} to deduce a primal-dual pair $(\bDelta,U_rV_r^T,\bU_{p_1 - r}\bV_{p_2 - r}^T)$ that satisfies
\begin{gather*}
\hSigma(\Thetas + \bDelta) - q + \lambda(U_rV^T + \bU_{p_1 - r}\bV_{p_2 - r}^T) = 0 \\
\bU_{p_1 - r}\bV_{p_2 - r}^T \in I
\end{gather*}
is unique, consistent, and PDW feasible.

\begin{lemma}
\label{lem:linearized-low-rank-multivariate-regression-consistency}
Assume the linearized problem \eqref{eq:nuclear-norm-minimization-linearized-cov} satisfies RSC (over $\linspan(I)^\perp$) and \eqref{eq:low-rank-multivariate-regression-irrepresentability}. For $\lambda = \frac{8(2-\alpha)}{\alpha}\sigma\sqrt{\frac{p_1+p_2}{n}},$ the unique primal-dual pair for \eqref{eq:nuclear-norm-minimization-linearized-cov} $(\bDelta,U_rV_r^T,\bU_{p_1 - r}\bV_{p_2 - r}^T)$ is
\BNUM
\item consistent: $\bigl\|\bDelta\bigr\|_{\Fro} \le \frac{4}{m}\left(\sqrt{2} + \frac{4(2-\alpha)}{\alpha}\right)\sigma\sqrt{\frac{r(p_1+p_2)}{n}}.$
\item PDW feasible: $\norm{\bU_{p_1 - r}\bV_{p_2 - r}^T}_{\op} \le 1 - \frac{\tau}{2}.$
\ENUM
\end{lemma}

By Theorem \ref{thm:dual-consistency} (and the convexity of the nuclear norm),
\begin{align*}
&\bigl\|\hU\hV^T - U_r V_r^T - \bU_{p_1 - r}\bV_{p_2 - r}^T\bigr\|_{\op}^2 \\
&\pc\le \|\hU\hV^T - U_r V_r^T - \bU_{p_1 - r}\bV_{p_2 - r}^T\|_{\Fro}^2 \\
&\pc\le \frac{M}{\lambda}\bigl(R(\bDelta) - R(\hDelta)\bigr) \le \frac{M}{\lambda}R(\bDelta),
\end{align*}
where $M = \sup_{\norm{\Delta}_{\Fro}\,\le\,1}\norm{Q\Delta}_{\Fro}.$ Since $\bDelta\in\linspan(I)^\perp,$ we may apply Lemma \ref{lem:nuclear-norm-taylors-thm} to obtain
\[
\bigl\|\hU\hV^T - U_r V_r^T - \bU_{p_1 - r}\bV_{p_2 - r}^T\bigr\|_{\op}^2 \le \frac{4}{3s_r}\frac{M}{\lambda}\|\bDelta\|_{\Fro}^2
\]
as long as $\|\bDelta\|_{\op} \le \frac{s_r}{2}.$ By the consistency of the linearized problem
\[
\|\bDelta\|_{\op} \le \|\bDelta\|_{\Fro} \le \frac{4}{m}(\sqrt{2} + \alpha')\sigma\sqrt{\frac{r(p_1+p_2)}{n}},
\]
where $\alpha' = \frac{4(2-\alpha)}{\alpha}.$ We put the pieces together to obtain
\BEQ
\begin{aligned}
&\bigl\|\hU\hV^T - U_r V_r^T - \bU_{p_1 - r}\bV_{p_2 - r}^T\bigr\|_{\op}^2 \\
&\pc\le \frac{32}{3s_r}\frac{M}{m^2}\frac{(\sqrt{2} + \alpha')^2}{\alpha'}\sigma r\sqrt{\frac{(p_1+p_2)}{n}},
\end{aligned}
\label{eq:nuclear-norm-minimization-rank-consistency-1}
\EEQ
when $n > \frac{16}{m^2}\frac{\sigma^2}{s_r^2}(\sqrt{2} + \alpha')^2r(p_1+p_2).$

By Lemma \ref{lem:linearized-problem-consistency}, $\bU_{p_1 - r}\bV_{p_2 - r}^T$ is PDW feasible. Thus it has at most $r$ unit singular values. Its $\min\{p_1,p_2\} - r$ remaining singular values are smaller than $1 - \frac{\alpha}{2}.$ By Weyl's inequality, it suffices to ensure
\BEQ
\bigl\|\hU\hV^T - U_r V_r^T - \bU_{p_1 - r}\bV_{p_2 - r}^T\bigr\|_{\op} \le \frac{\alpha}{2}.
\label{eq:nuclear-norm-minimization-rank-consistency-2}
\EEQ
to ensure $\bU\bV^T$ has no more than $r$ unit singular values. We combine \eqref{eq:nuclear-norm-minimization-rank-consistency-1} and \eqref{eq:nuclear-norm-minimization-rank-consistency-2} to deduce the requirement on $n.$
\end{proof}

To our knowledge, Theorem \ref{cor:nuclear-norm-minimization-rank-consistency} is the first non-asymptotic rank consis\-tency result for nuclear norm minimization. To keep things simple, we considered low-rank multivariate regression. However, the proof technique generalizes readily to M-estimators with other loss functions and regularizers.

\section{Computational experiments}

\label{sec:experiments}

We show some consequences of Corollary \ref{cor:mle-consistent} with experiments on two models from structure learning of networks that are motivated by bioinformatics applications. We select $\lambda$ proportional to $\sqrt{ (\max_{g\in\cG}|g|)\frac{\log|\cG|}{n}}$ and use a proximal Newton-type method \cite{lee2012proximal} to solve the likelihood maximization problem.

\subsection{The graphical lasso}

Suppose we are given samples drawn \iid\ from a normal distribution, and we seek to estimate the inverse covariance matrix. If $p > n$, we cannot simply invert the sample covariance matrix $\hat{\Sigma}$ to estimate $\Theta^\star$. However, we can obtain a maximum likelihood estimate of $\Theta^\star$:
\begin{align}
  \minimize_\Theta\;\tr\left(\hat{\Sigma}\Theta\right) - \log\det(\Theta) + \lambda\sum_{s,t\in\cG}\norm{\Theta_{st}}_2.
  \label{eq:l1-logdet}
\end{align}
The group lasso penalty to promotes block sparse inverse covariance matrices, and $\lambda$ trades-off goodness-of-fit and group sparsity. This estimator is a group sparse variant of the \emph{graphical lasso} \cite{friedman2008sparse}.



\begin{figure}
  \begin{subfigure}{0.5\textwidth}
  \includegraphics[width=\textwidth]{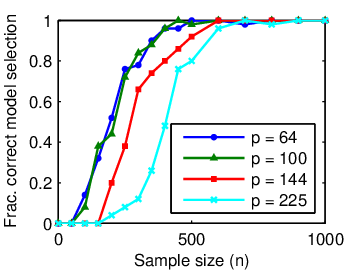}
  \end{subfigure}%
   ~
  \begin{subfigure}{0.5\textwidth}
   \includegraphics[width=\textwidth]{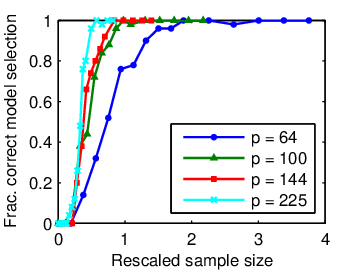}
  \end{subfigure}
  \caption{Fraction of correct model selection versus sample size $n$ and rescaled sample size $n/((\max_{g\in\cG}|g|)\log|\cG|)$ with the grouped graphical lasso. Each point represents the fraction of 100 trials when the grouped graphical lasso correctly estimated the true group structure.}
  \label{fig:group-glasso}
\end{figure}

We estimate the probability of correct model selection using the fraction of 100 trials when the graphical lasso correctly estimates the true group structure. Figure \ref{fig:group-glasso} shows the fraction of correct group structure selection versus the sample size $n$ for four graphs consisting of 64, 100, 144, and 225 nodes. In these experiments, we varied the sample size $n$ from 100 to 1000.

The fraction of correct model selection is small for small sample sizes but grows to one as the sample size increases. Intuitively, more samples are required to learn a larger model, hence the curves for larger graphs are to the right of curves for smaller graphs.
If we plot these curves with the x-axis rescaled by $1/((\max_{g\in\cG}|g|)\log|\cG|)$, then the curves align. This is consistent with Corollary \ref{cor:penalized-mle-consistent} that say the effective sample size scales logarithmically with $|\cG|$.

\subsection{Learning mixed graphical models}

The pairwise mixed graphical model was developed to model data that contain both categorical and continuous features \cite{lee2012learning, cheng2013high} \eg, two features about a person are weight (continuous) and gender (categorical). The model is a natural pairwise extension of the Gaussian MRF and a pairwise discrete MRF:
\BEQ
\Pr(x,y;(\beta,\theta,\gamma)) \propto \exp\biggl(\sum_{s,t}-\textstyle\frac12 \beta_{st} x_{s} x_{t}+\displaystyle\sum_{s,j} \theta_{sj}(y_{j})x_{s} +\sum_{j,r} \gamma_{rj}(y_r , y_j)\Big).
\EEQ
$x_s,\,s = 1,\dots,p$ and $y_j,\,j = 1,\dots,q$'s are continuous and discrete variables and $\beta_{st},\theta_{sj},\gamma_{rj}$ are continuous-continuous, continuous-discrete, and discrete-discrete edge potentials. We seek maximum likelihood and pseudolikelihood estimates of the parameters $(\beta,\theta,\gamma)$
\begin{align}
\minimize_{(\beta,\theta,\gamma)}\; &-\ell^{(n)}((\beta,\theta,\gamma)) + \lambda\rho((\beta,\theta,\gamma)).
\label{eq:mixed-lasso}
\end{align}
$\rho$ is the group lasso penalty:
$$
\rho((\beta,\theta,\gamma)) = \sum_{s,t} \abs{\beta_{st}}  +\sum_{s,j} \norm{\theta_{sj}}_2 +\sum_{j,r} \norm{\gamma_{rj}}_F.
$$
To make sure the model is identifiable, we enforce linear constraints on $\gamma_{rj}$:
\[\textstyle
\sum_{x_r, x_j} \gamma_{rj} (x_r,x_j) = 0,\,j,r = 1,\dots,q.
\]

\begin{figure}
  \begin{subfigure}{0.5\textwidth}
  \centering
  \includegraphics[width=0.75\textwidth]{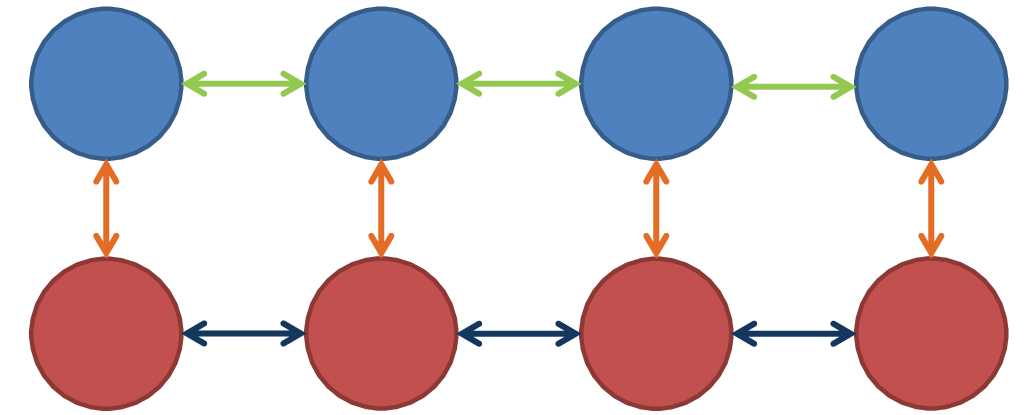}
  \subcaption{The graph topology used in this experiment. The blue nodes are continuous variables and the red nodes are discrete variables. The actual experiment had $10$ continuous and $10$ discrete variables.}
  \label{fig:syntheticgraph}
  \end{subfigure}%
   ~
  \begin{subfigure}{0.5\textwidth}
   \includegraphics[width=\textwidth]{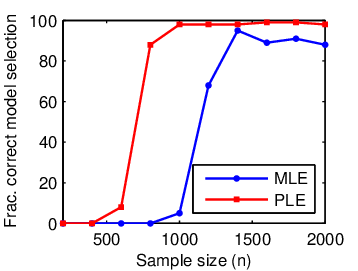}
  \end{subfigure}
  \caption{Fraction of correct model selection versus sample size $n$ of the penalized MLE and PLE on a mixed graphical model. Each point represents the fraction of 100 trials when the grouped graphical lasso correctly estimated the true group structure.}
  \label{fig:mixed-model}
\end{figure}

We estimate the probability of correct model selection using the fraction of 100 trials when \eqref{eq:mixed-lasso} correctly estimates the true group structure. Figure \ref{fig:mixed-model} shows the fraction of correct group structure selection versus the sample size $n$. In these experiments, we varied the sample size from 300 to 2000.

The fraction of correct model selection is small for small sample sizes but grows with the sample size. For the penalized PLE, the fraction grows to one, but, for the penalized MLE, the fraction plateaus around 0.9. This can be explained by the penalized MLE violating the irrepresentable condition. We refer to Section 3.1.1 in \cite{ravikumar2011high} for a similar example where the the irrepresentable condition holds for a neighborhood-selection estimator but fails for the penalized MLE.

\section{Conclusion}

We proposed the notion of geometric decomposablility and showed it is key to the model selection properties of regularized M-estimators. Our notion of decomposability builds on the notions by \cite{candes2012simple} and \cite{vandegeer2012weakly} and readily admits a notion of irrepresentability.

We also developed a general framework for establishing consistency and model selection consistency of regularized M-estimators. Our main result (Theorem \ref{thm:main} gives deterministic conditions on the problem that guarantee consistency and model selection consistency. We combined our main result with probabilistic analysis to study the model selection properties of the lasso, generalized lasso, and nuclear norm minimization. To our knowledge the non-asymptotic result on rank-consistency of nuclear norm minimization is novel.

\section*{Acknowledgements}

We thank Trevor Hastie and three anonymous reviewers for their insightful comments. J. Lee was supported by a National Defense Science and Engineering Graduate Fellowship (NDSEG), a NSF Graduate Fellowship, and a Stanford Graduate Fellowship. Y. Sun was supported by the NIH, grant 1U01GM102098-01. J.E. Taylor was supported by the NSF, grant DMS 1208857, and by the AFOSR, grant 113039.

\appendix

\section{Proof of lemmas in Section 3}
\label{sec:main-result-proofs}

\begin{proof}[Proof of Lemma \ref{lem:V-semi-norm}]
$V$ is the infimal convolution of $\gamma_I = h_{I^\circ}$ and $\ones_E = h_{E^\perp}.$ By the properties of support functions, $V = h_{E\,\cap\,I^\circ}.$ Since $I^\circ$ is bounded and support functions of bounded sets are finite and sublinear, $V$ is finite and sublinear.
\end{proof}

\begin{proof}[Proof of Lemma \ref{lem:consistency}]
Since $\btheta$ solve the restricted problem, we have
$$
\ell(\btheta) + \lambda h_A(\btheta) \le \ell(\theta^\star) + \lambda h_A(\theta^\star).
$$
Since $\hat \theta \in C$ and the objective is strongly convex over $C$, $\hat \theta$ is the unique solution to \eqref{eq:regularized-m-estimator-1}. By Assumption \ref{asu:rsc}, we have
\[
W^T P_M(\btheta - \theta^\star) +\frac{m}{2} \|\btheta-\theta^\star\|_2 ^2 + \lambda(\rho(\btheta) - \rho(\theta^\star)) \le 0,
\]
where $W = \nabla \ell(\theta^\star).$ We take norms to obtain
\begin{align}
0 &\ge - \varrho^*(P_MW)\varrho\bigl(\btheta - \theta^\star\bigr) + \frac{m}{2}\bigl\|\btheta - \theta^\star\bigr\|_2^2 - \lambda \rho(\btheta - \theta^\star) \\
&\ge - \kappa_{\varrho}\varrho^*(P_MW)\bigl\|\btheta - \theta^\star\bigr\|_2 + \frac{m}{2}\bigl\|\btheta - \theta^\star\bigr\|_2^2 - \lambda \rho(\btheta - \theta^\star).
\label{eq:consistency-2}
\end{align}
Since $\htheta - \thetas \in M,$ we have
$$
h_A(\btheta - \theta^\star) = \rho(\btheta - \theta^\star) \le \kappa_\rho\bigl\|\btheta - \theta^\star\bigr\|_2.
$$
We substitute this bound into \eqref{eq:consistency-2} to obtain
\begin{align*}
0 \ge - \kappa_{\varrho}\varrho^*(P_MW)\bigl\|\btheta - \theta^\star\bigr\|_2 + \frac{m}{2}\bigl\|\btheta - \theta^\star\bigr\|_2^2 - \kappa_\rho\lambda\bigl\|\btheta - \theta^\star\bigr\|_2.
\end{align*}
This means
$$
\bigl\|\btheta - \theta^\star\bigr\|_2 \le \frac{2}{m}\left(\kappa_{\varrho}\varrho^*(P_MW) + \kappa_\rho\lambda\right).
$$
Plugging in the choice of $\lambda > \frac{4\kappa_{\IC}}{\tau}\varrho^*(P_MW)$ gives our conclusion.
\end{proof}

\begin{proof}[Proof of Proposition \ref{lem:dual-certificate}]
Suppose there are two primal dual solution pairs, $(\theta_1,z_{A,1},z_{I,1},z_{E^\perp,1})$ and $(\theta_2,z_{A,2},z_{I,2},z_{E^\perp,2})$, \ie{}
\begin{gather}
\nabla \ell(\theta_1) + \lambda(z_{A,1} + z_{I,1} + z_{E^\perp,1}) = 0
\label{eq:dual-certificate-1} \\
\nabla \ell(\theta_2) + \lambda(z_{A,2} + z_{I,2} + z_{E^\perp,2}) = 0.
\end{gather}
Since the original problem \eqref{eq:regularized-m-estimator-1} is convex, the optimal value is unique:
\begin{align*}
&\ell(\theta_1) + P(\theta_1) = \ell(\theta_1) + \lambda(z_{A,1} + z_{I,1} + z_{E^\perp,1})^T\theta_1 \\
&\pc= \ell(\theta_2) + P(\theta_2) = \ell(\theta_2) + \lambda(z_{A,2} + z_{I,2} + z_{E^\perp,2})^T\theta_2.
\end{align*}
We subtract $\lambda(z_{A,1} + z_{I,1} + z_{E^\perp,1})^T\theta_2$ from both sides to obtain
\begin{align*}
&\ell(\theta_1) + \lambda(z_{A,1} + z_{I,1} + z_{E^\perp,1})^T(\theta_1 - \theta_2) \\
&\pc= \ell(\theta_2) + \lambda(z_{A,2} + z_{I,2} + z_{E^\perp,2} - z_{A,1} - z_{I,1} - z_{E^\perp,1})^T\theta_2.
\end{align*}
We rearrange this expression to obtain
\begin{align*}
&\ell(\theta_1) - \ell(\theta_2) + \lambda(z_{A,1} + z_{I,1} + z_{E^\perp,1})^T(\theta_1 - \theta_2)\\
&\pc= \lambda(z_{A,2} + z_{I,2} + z_{E^\perp,2} - z_{A,1} - z_{I,1} - z_{E^\perp,1})^T\theta_2.
\end{align*}
We substitute in \eqref{eq:dual-certificate-1} to obtain
\begin{align*}
&\ell(\theta_1) - \ell(\theta_2) - \nabla \ell(\theta_1)^T(\theta_1 - \theta_2)\\
&\pc= \lambda(z_{A,2} + z_{I,2} + z_{E^\perp,2} - z_{A,1} -z_{I,1} - z_{E^\perp,1})^T\theta_2.
\end{align*}
Since $\ell$ is convex, the left side is non-positive and
$$
(z_{A,2} + z_{I,2} + z_{E^\perp,2})^T\theta_2 \le (z_{A,1} + z_{I,1} + z_{E^\perp,1})^T\theta_2.
$$
Since $\theta_1$ and $\theta_2$ are in $S$, we can ignore the terms $z_{E^\perp,2}^T\theta_2$ and $z_{E^\perp,1}^T\theta_2$ to obtain
$$
(z_{A,2} + z_{I,2})^T\theta_2 \le (z_{A,1} +z_{I,1})^T\theta_2.
$$
But we also know
\begin{align*}
(z_{A,1} + z_{I,1})^T\theta_2 &\le \sup_u\,\{u^T\theta_2\mid u\in A\} + \sup_u\,\{u^T\theta_2\mid u\in I\} \\
&= z_{A,2}^T\theta_2 + z_{I,2}^T\theta_2.
\end{align*}
We combine these two inequalities to obtain
$$
(z_{A,2} + z_{I,2})^T\theta_2 = (z_{A,1} + z_{I,1})^T\theta_2 \le z_{A,2}^T\theta_2 + z_{I,1}^T\theta_2
$$
This simplifies to $z_{I,2}^T\theta_2 \le z_{I,1}^T\theta_2$. If $z_{I,1}\in\relint(I)$, then
\begin{align*}
z_{I,1}^T\theta_2 &= z_{I,2}^T\theta_2\,\textnormal{if }\theta_2\,\textnormal{has no component in }\linspan(I) \\
z_{I,1}^T\theta_2 &< z_{I,2}^T\theta_2\,\textnormal{if }\theta_2\,\textnormal{has a component in }\linspan(I).
\end{align*}
But we also know $z_{I,2}^T\theta_2 \le z_{I,1}^T\theta_2$. Thus $\theta_2$ has no component in $\linspan(I)$.
\end{proof}

\begin{proof}[Proof of Lemma \ref{lem:R-small}]
The Taylor remainder term is simply
$$
R = \nabla\ell(\btheta) - \nabla\ell(\theta^\star) - Q(\btheta - \theta^\star).
$$
By mean value theorem (along $\btheta - \thetas$), we have
$$
R = \int_0^1\bigl(\nabla^2\ell(\theta^\star + \alpha(\btheta - \theta^\star)) - Q\bigr)(\btheta - \theta^\star)\,d\alpha.
$$
Since $\nabla^2\ell$ is Lipschitz continuous with constant $L$ over $C$,
\begin{align*}
\|R\|_2 &= \norm{\int_0^1\bigl(\nabla^2\ell(\theta^\star + \alpha(\btheta - \theta^\star)) - Q\bigr)(\btheta - \theta^\star)\,d\alpha}_2 \\
&\le \int_0^1\bigl\|\nabla^2\ell(\theta^\star + \alpha(\btheta - \theta^\star)) - Q\bigr\|\bigl\|\btheta - \theta^\star\bigr\|_2\,d\alpha \\
&\le \int_0^1L \bigl\|\btheta - \theta^\star\bigr\|_2^2 \alpha \,d\alpha \\
&\le \frac{L}{2}\bigl\|\btheta - \theta^\star\bigr\|_2^2.
\end{align*}
By Lemma \ref{lem:consistency}, we have
\[
\|R\|_2 \le \frac{2L}{m^2}\Bigl(\kappa_\rho + \frac{\tau}{4}\frac{\kappa_\varrho}{\kappa_{\IC}}\Bigr)^2\lambda^2.
\]
To ensure $\frac{\kappa_{\IC}}{\lambda}\varrho^*(R) \le \frac{\tau}{4},$ it sufficies to ensure $\frac{\kappa_{\IC}}{\lambda}\norm{R}_2 \le \frac{\tau}{4\kappa_{\varrho^*}}.$ Plugging in the choice of $\lambda$ gives the desired conclusion.
\end{proof}

\section{Proofs of lemmas in Section 4}
\label{sec:examples-proofs}

\begin{proof}[Proof of Lemma \ref{lem:subexp-concentration}]

First, we prove an auxiliary result: the sufficient statistics of an exponential family are subexponential random variables.

\begin{lemma}
\label{lem:suff-stat-subexp}
Let $\phi^\star_i = \Expect_{ \theta^\star} \left[ \phi_i (x)\right]$ and $\phi_i (x)$ be a sufficient statistics of a regular exponential family. The random variable $\phi_i (x) - \phi^\star_i$ is subexponential :
$$
\Expect_{\theta^\star}\left[ \exp{s_i ( \phi_i (x) -\phi^\star_i) }\right]\le \exp(1)
$$
for some $s_i>0$.
\end{lemma}

\begin{proof}
\begin{align*}
&\Expect_{\theta^\star}\left[ \exp{s_i( \phi_i (x) -\phi^\star_i) }\right]\\
&\int dx \ h(x) \exp\left( \theta^{\star T} \phi(x) +s_i \phi_i (x) - A(\theta^\star) - s_i \phi^\star_i \right)\\
&\exp\left( -A(\theta^\star) -s_i \phi^\star_i \right) \int dx \ h(x) \exp \left( (\theta^\star +s_i e_i)^T \phi(x) \right)\\
&\exp\left( -A(\theta^\star) -s_i \phi^\star_i +A(\theta^\star +s_i e_i ) \right) \\
&\exp(-s_i \phi^\star_i) \exp\left( A(\theta^\star +s_i e_i )-A(\theta^\star)\right).
\end{align*}
Using continuity and the regular exponential family,  $|A(\theta^\star+s_i e_i)-A(\theta^\star)| <\epsilon$ and $|-s_i \phi^\star_i| <\epsilon$ for small enough $s_i$.  Thus
$$
 \exp\left( A(\theta^\star +s_i e_i )-A(\theta^\star)-s_i \phi^\star_i\right) \le\exp(1)
$$
for small enough $s_i$.
\end{proof}

We have  $\frac{\partial \ell}{\partial \theta_j} (\theta^\star) = \frac{1}{n}\sum_{i=1}^n \left( -\phi_j (x^{(i)}) + \Expect_{ \theta^\star} \left[ \phi_j (x)\right]\right)$. Thus $\frac{\partial \ell}{\partial \theta_j} (\theta^\star) $ is a sum of \iid\ subexponential random variables (Lemma \ref{lem:suff-stat-subexp}) and applying \cite[Corollary 5.17]{vershynin2010introduction} gives
$$
\Pr \left(\left|\frac{\partial \ell}{\partial \theta_j} (\theta^\star)\right| >t \right) \le 2\exp\left( -cn \min(t^2 /K_j ^2, t/K_j)\right)
$$
where $K_j$ is the Orlicz 1-norm of $-\phi_j (X) + \Expect_{ \theta^\star} \left[ \phi_j (X)\right]$ \cite[Definition 5.13]{vershynin2010introduction}. Let $K = \max_j K_j$. By the union bound,
\begin{align*}
\Pr \left(\norm{\nabla_{\theta_g} \ell (\theta^\star) }_2 >t \right) &\le \Pr \left( \text{for some } j \in g,\ \left|\frac{\partial \ell}{\partial \theta_j} (\theta^\star)\right| >t/\sqrt{|g|}\right)\\
&\le \sum_{j\in g} \Pr \left(\left|\frac{\partial \ell}{\partial \theta_j} (\theta^\star)\right| >t/\sqrt{|g|} \right)\\
&\le 2\exp\left( -cn \min\left(\frac{t^2} { |g| K ^2}, \frac{t}{ \sqrt{|g|}K}\right)\right)\\
&\le 2\exp\left( -cn \min\left(\frac{t^2} { \max_{g \in \cG} |g| K ^2}, \frac{t}{ \max_{g \in \cG} \sqrt{|g|}K}\right)\right)
\end{align*}
and
\begin{align*}
&\Pr \left(\max_{g \in \cG} \norm{\nabla_{\theta_g} \ell (\theta^\star) }_2 >t \right) \le \Pr \left( \text{for some } g \in \cG,\ \norm{\nabla_{\theta_g} \ell (\theta^\star) }_2 >t \right) \\
&\le \sum_{g \in \cG} \Pr \left(\norm{\nabla_{\theta_g} \ell (\theta^\star) }_2 >t \right)\\
&\le |\cG| 2\exp\left( -cn \min\left(\frac{t^2} { \max_{g \in \cG} |g| K ^2}, \frac{t}{ \max_{g \in \cG} \sqrt{|g|}K}\right)\right)\\
&=2\exp\left(\log |\cG| -cn \min\left(\frac{t^2} { \max_{g \in \cG} |g| K ^2}, \frac{t}{ \max_{g \in \cG} \sqrt{|g|}K}\right)\right)
\end{align*}
\end{proof}

\begin{proof}[Proof of Lemma \ref{lem:exp-family-compact-subset}]
By optimality of $\hat \theta$,
\begin{align*}
\ell(\hat \theta) +\lambda \bigl\|\hat{\theta}\bigr\|_{2/1} &\le \ell( \theta^\star) +\lambda \norm{ \theta^\star}_{2/1}\\
- \hat \theta^T \phi^n + A(\hat \theta)+\lambda \bigl\|\hat{\theta}\bigr\|_{2/1} &\le - \theta^{\star T}\phi^n + A( \theta^{\star})+\lambda \norm{ \theta^\star}_{2/1}\\
- \hat \theta^T \phi^n +\nabla A(\theta^\star)^T (\hat\theta - \theta^\star)+\lambda \bigl\|\hat{\theta}\bigr\|_{2/1} &\le - \theta^{\star T}\phi^n +\lambda \norm{ \theta^\star}_{2/1}\\
- \hat \theta^T \phi^n + \phi^{\star T} (\hat\theta - \theta^\star) +\lambda \bigl\|\hat{\theta}\bigr\|_{2/1} &\le - \theta^{\star T}\phi^n +\lambda \norm{ \theta^\star}_{2/1}\\
\lambda \bigl\|\hat{\theta}\bigr\|_{2/1} &\le \lambda \norm{ \theta^\star}_{2/1}+ (\hat\theta -\theta^{\star})^T\left( \phi^n - \phi^\star \right)
\end{align*}
We now bound $\bigl\|\hat{\theta}\bigr\|_{2/1}$,
\begin{align*}
\lambda \bigl\|\hat{\theta}\bigr\|_{2/1} &\le \lambda \norm{ \theta^\star}_{2/1}+ (\hat\theta -\theta^{\star})^T\left( \phi^n - \phi^\star \right)\\
&\le\lambda \norm{ \theta^\star}_{2/1}+\bigl\|\hat\theta -\theta^{\star}\bigr\|_{2/1}\norm{ \phi^n - \phi^\star}_{2,\infty}\\
&\le\lambda \norm{ \theta^\star}_{2/1}+\bigl\|\hat{\theta}\bigr\|_{2/1}\norm{ \phi^n - \phi^\star}_{2,\infty} + \norm{\theta^{\star}}_{2/1}\norm{ \phi^n - \phi^\star}_{2,\infty}
\end{align*}
Rearranging gives us,
\begin{align}
 \bigl\|\hat{\theta}\bigr\|_{2/1}  \le\frac{1}{(\lambda-\norm{ \phi^n - \phi^\star}_{2,\infty}) }\left(   \lambda \norm{ \theta^\star}_{2/1}+ \norm{ \phi^n - \phi^\star}_{2,\infty}\norm{\theta^{\star}}_{2/1}\right)
 \label{eq:compact-subset-bound}
\end{align}
For the second part,
\begin{align*}
\ell(\hat \theta) +\lambda \bigl\|\hat{\theta}\bigr\|_{2/1} &\le \ell( \theta^\star) +\lambda \norm{ \theta^\star}_{2/1}\\
- \hat \theta^T \phi^n + A(\hat \theta)+\lambda \bigl\|\hat{\theta}\bigr\|_{2/1} &\le - \theta^{\star T}\phi^n + A( \theta^{\star})+\lambda \norm{ \theta^\star}_{2/1}\\
A(\hat \theta) &\le (\hat \theta - \theta^\star)^T \phi^n +A(\theta^\star) +\lambda \norm{\theta^\star}_{2/1} -\lambda \bigl\|\hat{\theta}\bigr\|_{2/1}\\
A(\hat \theta) &\le \bigl\|\hat{\theta}\bigr\|_{2/1} \norm{\phi^n}_{2,\infty}+\norm{\theta^\star}_{2/1} \norm{\phi^n}_{2,\infty} +A(\theta^\star) +\lambda\norm{\theta^\star}_{2/1}.
\end{align*}
\end{proof}

\begin{proof}[Proof of Corollary \ref{cor:regularized-mle-consistent}]
Let $f(n,|\cG|,|g|)$ be a function that inverts the concentration inequality of Lemma \ref{lem:subexp-concentration} in the sense
\begin{align*}
&\Pr \left(\max_{g \in \cG} \norm{\nabla_{\theta_g} \ell (\theta^\star) }_2 >f(n,|\cG|,|g|) \right)\\
&\le 2\exp\left(\log |\cG| -cn \min\left(\frac{f(n,|\cG|,|g|)^2} { \max_{g \in \cG} |g| K ^2}, \frac{f(n,|\cG|,|g|)}{ \max_{g \in \cG} \sqrt{|g|}K}\right)\right)\\
&= 2 \exp(0).
\end{align*}
Thus $f$ is chosen so
\begin{align}
\log |\cG| - cn \min\left( \left(\frac{f}{\max_{g \in G} \sqrt{|g|} K}\right)^2, \frac{f}{\max_{g \in G} \sqrt{|g|} K}\right)=0
\label{eq:invert-concentration}
\end{align}
Let
$$
f(n,|\cG|,|g|) := \max_{g \in G} \sqrt{|g|} K \sqrt{\frac{\log |\cG|}{cn}}.
$$
For $n>\left(\frac{3}{2}\right)^2 \frac{\log |\cG|}{c}$ the first term in the min is active, so the choice of $f$ satisfies \eqref{eq:invert-concentration}.

By the following computation, the choice $\lambda = \frac{3 \kappa_{\IC}}{\tau}f(n,|\cG|,|g|)$ ensures that $$
\Pr \left(  \frac{2 \kappa_{\IC}}{\tau} \max_{g \in \cG} \norm{\nabla_{\theta_g} \ell (\theta^\star) }_2 >\lambda \right)< 2 |\cG|^{-5/4}
$$

\begin{align*}
&\Pr \left(  \frac{2 \kappa_{\IC}}{\tau} \max_{g \in \cG} \norm{\nabla_{\theta_g} \ell (\theta^\star) }_2 >\lambda \right)\\
&=\Pr \left(  \frac{2 \kappa_{\IC}}{\tau} \max_{g \in \cG} \norm{\nabla_{\theta_g} \ell (\theta^\star) }_2 >\frac{3 \kappa_{\IC}}{\tau} \max_{g \in G} \sqrt{|g|} K \sqrt{\frac{\log |\cG|}{cn}}\right)\\
&=\Pr \left(   \max_{g \in \cG} \norm{\nabla_{\theta_g} \ell (\theta^\star) }_2 >\frac{3 }{2} \max_{g \in G} \sqrt{|g|} K \sqrt{\frac{\log |\cG|}{cn}}\right)\\
&\le 2 \exp \left( \log |\cG| - cn\min\left( 9/4 \max_{g \in G} |g| K^2 \frac{\log |\cG|}{cn}, 3/2 \frac{f}{\max_{g \in G} \sqrt{|g|} K}\right)\right)\\
&= 2 \exp \left( \log |\cG| - cn\min\left( 9/4  \frac{\log |\cG|}{cn}, 3/2 \sqrt{\frac{\log |\cG|}{cn}}\right)\right)
\end{align*}
Since $n>\left(\frac{3}{2}\right)^2 \frac{\log |\cG|}{c}$, we have $9/4  \frac{\log |\cG|}{cn}< 3/2 \frac{\log |\cG|}{cn}$ and thus
\begin{align*}
&\Pr \left(  \frac{2 \kappa_{\IC}}{\tau} \max_{g \in \cG} \norm{\nabla_{\theta_g} \ell (\theta^\star) }_2 >\lambda \right)\\
&\le2 \exp \left( \log |\cG| - \frac{9}{4}\log |\cG|\right)\\
& = 2 |\cG|^{-5/4}.
\end{align*}

For the rest of this proof, we will assume the event $\{\lambda > \frac{2 \kappa_{\IC}}{\tau}\max_{g \in \cG} \norm{\nabla_{\theta_g} \ell (\theta^\star) }_2\}$, so all the following statements hold with probability at least $1- 2 |\cG|^{-5/4}$.

Lemma \ref{lem:exp-family-compact-subset} shows
\begin{align*}
 &\bigl\|\hat{\theta}\bigr\|_{2/1}  \le\frac{   \lambda \norm{ \theta^\star}_{2/1}+ \norm{ \phi^n - \phi^\star}_{2,\infty}\norm{\theta^{\star}}_{2/1}}{\lambda-\norm{ \phi^n - \phi^\star}_{2,\infty}}\\
 &\le \frac{   \lambda \norm{ \theta^\star}_{2/1}+ \frac{\tau}{2\kappa_{\IC}} \lambda \norm{\theta^{\star}}_{2/1}}{\lambda-\frac{\tau}{2\kappa_{\IC}} \lambda}\\
 &\le \frac{2\norm{\theta^{\star}}_{2/1}} { 1- \frac{\tau}{2\kappa_{\IC}}}\\
 &\le 4 \norm{\theta^{\star}}_{2/1}
\end{align*}
where we used $\norm{ \phi^n - \phi^\star}_{2,\infty} = \max_{g \in \cG} \norm{\nabla_{\theta_g} \ell (\theta^\star) }_2$ and $\frac{\tau}{\kappa_{\IC}}\le 1$.
Lemma \ref{lem:exp-family-compact-subset} also shows that
\begin{align*}
 A(\hat \theta) &\le \norm{\theta^{\star}}_{2/1} \norm{\phi^n}_{2,\infty} + \bigl\|\hat{\theta}\bigr\|_{2/1} \norm{\phi^n}_{2,\infty} +A(\theta^\star) +\lambda \norm{\theta^\star}_{2/1}\\
 &\le  5\norm{\theta^{\star}}_{2/1} \norm{\phi^n}_{2,\infty} +A(\theta^\star) +3\frac{\kappa_{\IC}}{\tau} f(n,|\cG|,|g|) \norm{\theta^\star}_{2/1}\\
 &\le 5\norm{\theta^{\star}}_{2/1} \left(\norm{\phi^\star}_{2,\infty} +\frac{3}{2}f(n,|\cG|,|g|)\right) +A(\theta^\star) +3\frac{\kappa_{\IC}}{\tau} f(n,|\cG|,|g|) \norm{\theta^\star}_{2/1}\\
 &=5\norm{\theta^{\star}}_{2/1} \left(\norm{\phi^\star}_{2,\infty} +\frac{3}{2}\max_{g \in G} \sqrt{|g|} K \sqrt{\frac{\log |\cG|}{cn}}\right) +A(\theta^\star) \\
 &+3\frac{\kappa_{\IC}}{\tau} \max_{g \in G} \sqrt{|g|} K \sqrt{\frac{\log |\cG|}{cn}} \norm{\theta^\star}_{2/1}\\
 &< =5\norm{\theta^{\star}}_{2/1} \left(\norm{\phi^\star}_{2,\infty} +\max_{g \in G} \sqrt{|g|} K\right) +A(\theta^\star) +2\frac{\kappa_{\IC}}{\tau} \max_{g \in G} \sqrt{|g|} K \norm{\theta^\star}_{2/1}\\
 &=:R
\end{align*}
by the triangle inequality, $\frac{3}{2}f(n,|\cG|,|g|) >\max_{g \in \cG}$, and $n>\left(\frac{3}{2}\right)^2 \frac{\log |\cG|}{c}$.

Thus from the above arguments we know that
$$
\hat \theta \in C:=\{\theta\ \big| \norm{\theta}_{2/1} \le 4 \norm{\theta^\star}_{2/1} \text{ and }A(\theta)\le R\}.
$$
The subset $C$ is convex and compact. Since the exponential family is minimal on $M$, $v^T \nabla^2 A(\theta) v >0$ for $v\in M$ \cite{wainwright2008graphical} and thus strongly convex over the compact subset $C \,\cap\, M$ with strong convexity constant $m$ (Assumption \ref{asu:rsc}). By the extreme value theorem applied to $\frac{\norm{\nabla^2 A(\theta) -\nabla^2 A(\theta^\star)}_2}{\norm{\theta -\theta^\star}_2}$, the function $\nabla^2\ell(\theta)$ has a finite Lipschitz constant $L$ over $C$.

Before we apply Theorem \ref{thm:main}, we compute the constants $\kappa_{\rho^*}, \kappa_{\rho^*}.$ Since the regularizer is finite (it's a norm), its dual semi-norm is finite. To keep things simple, we let $\varrho = \norm{\cdot}_{2/1}.$ The constants $\kappa_\rho = \kappa_\varrho, \kappa_{\rho^*}$ are
\begin{gather*}
\textstyle\kappa_\rho = \sup_\theta\left\{\norm{\theta}_{2/1}\mid \theta\in B_2\,\cap\,\linspan(B_{2/\infty,\cS})\right\} = \sqrt{\abs{\cS}}, \\
\textstyle\kappa_{\rho^*} = \sup_x\,\left\{\max_{g\in\cG}\mid\theta\in B_2\,\cap\,\linspan(B_{2/\infty,\cS})\right\} \le 1.
\end{gather*}

To apply Theorem \ref{thm:main}, we need to verify that the choice of $\lambda$ satisfies
$$
\lambda <\frac{m^2}{2L}\frac{\tau}{\kappa_{\IC}\kappa_{\rho^*}}\left(2\kappa_\rho + \frac{\kappa_\rho}{\kappa_{\IC}}\frac{\tau}{2}\right)^{-2}\frac{\tau}{2\kappa_{\IC}}.
$$
Substituting the expressions for the compatibility constants into the expression above gives
$$
\lambda <\frac{m^2}{2L}\left(2\sqrt{|\cS|} + \frac{\tau}{2\kappa_{\IC}}\sqrt{|\cS|}\right)^{-2}\frac{\tau}{2\kappa_{\IC}}
$$
or equivalently
\BEQ
\frac{3 \kappa_{\IC}}{\tau}f(n,|\cG|,|g|)<\frac{m^2}{L}\left(2\sqrt{|\cS|} + \frac{\tau}{2\kappa_{\IC}}\sqrt{|\cS|}\right)^{-2}\frac{\tau}{2\kappa_{\IC}} .
\label{eq:f-n-relation}
\EEQ

Using $ f(n,|\cG|,|g|)=\max_{g \in G} \sqrt{|g|} K \sqrt{\frac{\log |\cG|}{cn}}$,
\begin{gather*}
\frac{3 \kappa_{\IC}}{\tau} \max_{g \in G} \sqrt{|g|} K \sqrt{\frac{\log |\cG|}{cn}}<\frac{m^2}{L}\left(2\sqrt{|\cS|} + \frac{\tau}{2\kappa_{\IC}}\sqrt{|\cS|}\right)^{-2}\frac{\tau}{2\kappa_{\IC}} \\
\sqrt{cn} >6 \frac{ \kappa_{\IC}^2}{\tau^2}  \max_{g \in G} \sqrt{|g|} K \sqrt{\log |\cG|} \frac{L}{m^2} \left(2\sqrt{|\cS|} + \frac{\tau}{2\kappa_{\IC}}\sqrt{|\cS|}\right)^2\\
n> \frac{36}{c}\frac{ \kappa_{\IC}^4}{\tau^4}\max_{g \in G} |g| K^2\log |\cG|\frac{L^2}{m^4}\left(2\sqrt{|\cS|} + \frac{\tau}{2\kappa_{\IC}}\sqrt{|\cS|}\right)^4.
\end{gather*}

This completes the proof. We have verified all the assumptions of Theorem \ref{thm:main} and applying the theorem for the chosen value of $\lambda$ gives the desired result.
\end{proof}

\section{Proof of lemmas in Section 5}
\label{sec:weakly-decomposable-proofs}

\begin{proof}[Proof of Lemma \ref{lem:linearized-low-rank-multivariate-regression-consistency}]
Before we apply Lemma \ref{lem:linearized-problem-consistency}, we compute the constants $\kappa_{\brho},\kappa_{\bvarrho},\kappa_{\IC}.$ Since the regularizer is not a norm, we let $\bvarrho = \norm{\cdot}_*$ and check $\lambda > \frac{4\kappa_{\IC}}{\tau}\|\nabla\ell(\Thetas)\|_{\op}.$ It's straighforward to check
\[
\tr(U_r^T\Delta V) + \bigl\|U_{p_1-r}^T\Delta V_{p_2-r}\bigr\|_* \le \|\Delta\|_*.
\]
The ``model subspace'' $M$ is given by
\[
\linspan(I)^\perp = \left\{U_rX + (V_rY)^T\mid\text{ for any }X\in\reals^{r\times p_2},Y\in\reals^{r\times p_1}\right\},
\]
and the constants $\kappa_{\brho},\kappa_{\bvarrho}$ are given by
\begin{gather*}
\textstyle\kappa_{\brho} = \sup_{X,Y}\left\{\tr\left(U_r^T\Delta V_r\right)\mid \norm{\Delta}_{\Fro} \le 1\right\} = \sqrt{r} \\
\textstyle\kappa_{\bvarrho} = \sup_{X,Y}\left\{\norm{XV_r + U_rY^T}_*\mid\norm{U_rX + (V_rY)^T}_{\Fro} \le 1\right\} = \sqrt{2r}.
\end{gather*}
Similarly, the constant $\kappa_{\IC}$ is given by
\BEQ
\begin{aligned}
&\norm{U_{p_1-r}^T\bigl[P_I\bigl(QP_{I^\perp}(P_{I^\perp}QP_{I^\perp})^\dagger Z - Z\bigr)\bigr]V_{p_2-r}}_{\op} \\
&\pc\le \bigl\|U_{p_1-r}^T\bigl[P_IQP_{I^\perp}(P_{I^\perp}QP_{I^\perp})^\dagger Z\bigr]V_{p_2-r}\bigr\|_{\op} + \bigl\|U_{p_1-r}^TZV_{p_2-r}\bigr\|_{\op} \\
&\pc\le (2-\alpha)\norm{Z}_{\op}
\end{aligned}
\EEQ
is at most $2 - \alpha.$

To apply Lemma \ref{lem:linearized-problem-consistency}, we check $\lambda = \frac{8(2-\alpha)}{\alpha}\sigma\sqrt{\frac{p_1+p_2}{n}}$ satisfies the assumptions. By \cite{negahban2011estimation}, Lemma 3,
\[
\Pr\biggl(\frac{8(2-\alpha)}{\alpha n}\norm{\cX^*(\epsilon)}_2 > \frac{8(2-\alpha)}{\alpha}\sigma\sqrt{\frac{p_1 + p_2}{n}}\,\biggr) \le c_1e^{-c_2(p_1+p_2)},
\]
for some universal constants $c_1,c_2.$ Thus the claims of Lemma \ref{lem:linearized-problem-consistency} are valid with probability at least $1 - c_1e^{-c_2(p_1+p_2)}:$
\BNUM
\item consistent: $\bigl\|\bDelta\bigr\|_2 \le \frac{2}{m}\left(\sqrt{2} + \frac{4(2-\alpha)}{\alpha}\right)\sigma\sqrt{\frac{r(p_1+p_2)}{n}}.$
\item PDW feasible: $\norm{\bU_{p_1 - r}\bV_{p_2 - r}^T}_2 \le 1 - \frac{\tau}{2}.$
\ENUM
\end{proof}

\begin{proof}[Proof of Lemma \ref{lem:nuclear-norm-taylors-thm}]
For any $\Delta\in\linspan(I)^\perp,$ we have
\begin{align*}
&\norm{\Thetas + \Delta}_* - \norm{\Thetas}_* - \tr\bigl(V_rU_r^T\Delta \bigr) \\
&\pc= \tr\bigl(\tV_r\tU_r^T(\Thetas + \Delta) \bigr) - \tr\bigl(V_rU_r^T\Thetas \bigr) - \tr\bigl(V_rU_r^T\Delta \bigr),
\end{align*}
where $\tU\in\reals^{p_1\times r}$ and $\tV\in\reals^{p_2\times r}$ are the left and right singular factors of $\Thetas + \Delta.$ Since $\tr\bigl(\tV_r\tU_r^T\Thetas \bigr) \le \tr\bigl(V_rU_r^T\Thetas \bigr),$
\begin{align*}
&\norm{\Thetas + \Delta}_* - \norm{\Thetas}_* - \tr\bigl(V_rU_r^T\Delta \bigr) \le \tr\bigl(\bigl(\tU_r\tV_r^T - U_rV_r^T\bigr)^T\Delta \bigr) \\
&\pc\le \bigl\|\tU_r\tV_r^T - U_rV_r^T\bigr\|_{\Fro}\norm{\Delta}_{\Fro}.
\end{align*}
By \cite{li2002perturbation}, Theorem 2.4,
\[
\bigl\|\tU_r\tV_r^T - U_rV_r^T\bigr\|_{\Fro} \le \frac{4}{3s_r}\norm{\Delta}_F
\]
for any $\Delta$ such that $\norm{\Delta}_2 \le \frac{s_r}{2}.$ We put the pieces together to obtain the desired bound.

\end{proof}


\frenchspacing
\bibliographystyle{imsart-nameyear}
\bibliography{recovery}

\end{document}